\documentclass[a4paper,reqno,12pt]{amsart}
\usepackage[english]{babel}
\usepackage[T1]{fontenc}
\usepackage{newpxtext} 
\usepackage[a4paper,margin=3cm]{geometry}
\usepackage{amsmath,amsthm,amssymb,amsxtra,calligra,mathrsfs}
\usepackage{graphicx}
\usepackage{tikz}
\usepackage{tikz-cd} 
\usepackage{xcolor}
\usepackage{upgreek}
\usepackage{mathtools}
\usepackage{extarrows}
\usepackage{faktor}
\usepackage{tensor}
\usepackage{comment}
\usepackage{enumitem}
\usepackage{geometry}
\usepackage[all,cmtip]{xy}
\usepackage[colorlinks=true, backref=page]{hyperref}

\newcommand{\bZ}{\mathbb{Z}}
\newcommand{\bP}{\mathbb{P}}
\newcommand{\bC}{\mathbb{C}}

\newcommand{\bR}{\mathbb{R}}

\newcommand{\bD}{\mathbb{D}}

\newcommand{\bF}{\mathbb{F}}
\newcommand{\bQ}{{\mathbb Q}}
\newcommand{\bB}{{\mathbb B}}
\newcommand{\bH}{{\mathbb H}}

\newcommand{\cO}{{\mathcal O}}

\newcommand{\cY}{{\mathcal Y}}

\newcommand{\fm}{{\mathfrak m}}
\newcommand{\fM}{{\mathfrak M}}

\newcommand{\cA}{{\mathcal A}}

\newcommand{\sheafHom}{\mathscr{H}\text{\kern -3pt {\calligra\large om}}\,}

\DeclareMathOperator{\rk}{rk}

\DeclareMathOperator{\homo}{Hom}

\DeclareMathOperator{\Gr}{\mathrm{Gr}}

\DeclareMathOperator{\id}{\mathrm{id}}

\DeclareMathOperator{\MCG}{MCG}
\DeclareMathOperator{\PGL}{\mathrm{PGL}}
\DeclareMathOperator{\GL}{\mathrm{GL}}
\DeclareMathOperator{\SO}{\mathrm{SO}}
\DeclareMathOperator{\Sp}{\mathrm{Sp}}

\DeclareMathOperator{\NS}{\mathrm{NS}}
\DeclareMathOperator{\vol}{\mathrm{vol}}

\DeclareMathOperator{\aut}{Aut}

\DeclareMathOperator{\Alb}{Alb}

\newcommand{\q}{/\!\!/}

\newcommand\diag{\text{\rm diag}}

\newcommand\im{\text{\rm Im}}

\newcommand\rank{{\text{\rm rank}}}

\newcommand{\cF}{{\mathcal F}}

\newcommand\cD{{\mathcal{D}}}

\newcommand\cE{{\mathcal{E}}}

\newcommand\cC{{\mathcal{C}}}
\newcommand\cM{{\mathcal{M}}}

\newcommand\fX{{\mathfrak{X}}}

\newtheorem{defn-pro}{Definition-Proposition}
\newtheorem{defn-thm}{Definition-Theorem}
\newtheorem{thm}{Theorem}[section]
\newtheorem{lem}[thm]{Lemma}
\newtheorem{cor}[thm]{Corollary}
\newtheorem{pro}[thm]{Proposition}

\newtheorem{defn}[thm]{Definition}

\newtheorem{question}[thm]{Question}
\newtheorem{exmp}[thm]{Example}

\newtheorem{rem}[thm]{Remark}

\theoremstyle{remark}
\newtheorem*{remark}{Remark}

\title{On the global rigidity of period maps for certain families of varieties}
\author{ Cong Ding, Fei Si and Ruiran Sun }

\address{Cong Ding, School of Mathematical Sciences, Shenzhen University, 518060, Guangdong, China}
\email{congding@szu.edu.cn}

\address{Fei Si, The School of Mathematics and Statistics, Xi’an Jiaotong University, 28 West Xianning Road, Xi’an, Shaanxi 710049, China}
\email{sifei@xjtu.edu.cn}

\address{Ruiran Sun, School of Mathematical Sciences, Xiamen University, Xiamen 361005, China}
\email{ruiransun@xmu.edu.cn}

\date{}

\begin{document}

 \begin{abstract}
Motivated by Farb's question on the global rigidity of period maps and the effective Geometric Shafarevich Conjecture, we introduce the notion of \emph{$h$-deformation rigidity}: a holomorphic map is $h$-deformation rigid if it is the unique holomorphic map in its homotopy class. Using harmonic mapping theory, we develop a geometric criterion for this property and prove $h$-deformation rigidity for nonconstant holomorphic maps between equal-rank locally symmetric varieties of rank at least two   with irreducible source lattice  , as well as for certain ball-quotient inclusions and intermediate-Jacobian period maps.

 Then we study the rigidity problem for the period map via Kummer construction. We construct examples to show that the Kummer period map $\mathcal A_2 \to \mathcal F_4$ is not globally rigid. On the other hand, combining $h$-deformation rigidity with arithmetic monodromy and lattice-theoretic counting, we obtain an explicit estimate $\#\mathrm{Hol}_{\mathrm{nc}}(\mathcal A_2,\mathcal F_4)<2^{20249}$.
\end{abstract}
  \maketitle
  \setcounter{tocdepth}{1}
\tableofcontents
\section{Introduction }

In \cite{Farb22}, Farb proved a striking global rigidity theorem for the classical period map. More precisely, if $g>2$ and
\[
  F:\mathcal M_g\longrightarrow \mathcal A_h, \qquad h\leq g,
\]
is a nonconstant holomorphic map of complex orbifolds, then $h=g$ and $F$ is the Torelli period map. In his survey \cite{Farb24}, Farb proposed that one should study analogous rigidity phenomena for other period maps arising naturally in algebraic geometry.  This is the first motivation of the present paper. Let us formulate the problem more generally: suppose $M$ is a moduli space of certain varieties and there is a period map
\[
  p:M\longrightarrow \mathcal D/\Gamma
\]
to an arithmetic quotient of a bounded symmetric domain.  We assume throughout that the moduli spaces under consideration have orbifold singularities at worst, and all maps respect complex orbifold structure.   By a locally symmetric variety we mean an arithmetic quotient of a bounded symmetric domain.    We say that $p$ is \emph{globally rigid} if it is the unique nonconstant holomorphic orbifold map from $M$ to $\mathcal D/\Gamma$.  Thus Farb's theorem says that the Torelli map is globally rigid. See \cite{servan25, Zhong26} for other related global rigidity results.  The orbifold automorphisms $\aut_{orb}(M)$ and  $\aut_{orb}(\cD/\Gamma)$ provide a natural obstruction, so in some examples one should instead consider uniqueness modulo these automorphisms.

There is a second, somewhat different, motivation coming from the Geometric Shafarevich Conjecture.  The classical theorem of Parshin and Arakelov asserts the finiteness of non-isotrivial families of curves of fixed genus over a fixed curve with prescribed degeneracy \cite{Par68,Ara71}.  From this point of view, one can study not only the uniqueness of a distinguished period map, but also the \emph{cardinality} of moduli spaces of maps such as $\mathrm{Hol}_{\mathrm{nc}}(X,Y)$, where $X$ and $Y$ are moduli spaces or locally symmetric varieties. Moreover, we have particular interest on an \emph{explicit upper bound}.

The notion which relates these two problems in this paper is what we call \emph{homotopy-deformation rigidity}, or \emph{$h$-deformation rigidity} for short.  A holomorphic map $f:X\to Y$ is $h$-deformation rigid if every holomorphic map $g:X\to Y$ which is homotopic to $f$ is equal to $f$.  Thus an $h$-deformation rigid map is the unique holomorphic representative in its whole homotopy class. In particular, $h$-deformation rigidity implies the usual deformation rigidity.  Our first main result establishes this property in several general situations, where several natural examples of period maps arise in complex algebraic geometry could fit in.

\begin{thm}(=Theorems \ref{h-deform-rigidity_thmER}, \ref{h-deform-rigidity_thm_ball_case} and \ref{h-deform-rigidity_thm_intJac_case})
The following maps are $h$-deformation rigid:
\begin{enumerate}
    \item Any nonconstant holomorphic map $f:X\to Y$ between locally symmetric varieties with $ \operatorname{rk}(X)=\operatorname{rk}(Y)\geq 2$   and irreducible source lattice. Examples include the Kummer period map
    \[
      \tau_K:\mathcal A_2\longrightarrow \mathcal F_4.
    \]
    \item A holomorphic map $f$ from a ball quotient of dimension at least $2$ to the orthogonal locally symmetric variety $\mathcal D_{T}/\Gamma_0$ considered in Theorem~\ref{h-deform-rigidity_thm_ball_case} with $\mathrm{rk}(df) \ge 2$.
    \item The period map induced by intermediate-Jacobian
    \[
      M_{n,d}\longrightarrow \mathcal A_{N(n,d)}
    \]
    for the hypersurfaces listed in \S\ref{period:hypersurface}.
\end{enumerate}
\end{thm}

This theorem provides the first step towards both kinds of questions above.  For global rigidity, it separates the uniqueness problem inside a fixed homotopy class from the problem of determining which homotopy classes can occur.  For effective Shafarevich-type problems, it has an additional advantage: once $h$-deformation rigidity is known, counting holomorphic maps can be reduced to counting topological or arithmetic monodromy data.

Let us explain this point more precisely.  When the target is an orbifold $K(\pi,1)$, a free homotopy class is controlled by the conjugacy class of the induced homomorphism on orbifold fundamental groups.  Consequently, the global rigidity problem naturally separates into two problems:
\begin{enumerate}
    \item \emph{$h$-deformation rigidity}: prove uniqueness of a holomorphic map inside a fixed homotopy class;
    \item \emph{group rigidity}: classify the relevant homomorphisms
    \[
      \pi_1^{\mathrm{orb}}(M)\longrightarrow \Gamma
    \]
    up to conjugation in $\Gamma$.
\end{enumerate}
This proof strategy already appeared in Farb's proof of the global rigidity of the Torelli map \cite{Farb22}.  From this, we remark that $h$-deformation rigidity is a useful intermediate notion between deformation rigidity and global rigidity.

The notion of $h$-deformation rigidity is also inspired by the harmonic mapping approach to Shafarevich conjecture over complex function fields.  In particular, Jost and Yau used harmonic mappings to study algebraic varieties over function fields and Shafarevich-type conjecture in \cite{JostYau93}.  Farb's proof also uses the geodesic homotopy, together with an energy/Wirtinger argument on curves, to transform homotopy information into holomorphic rigidity.  Although the techniques in our paper    are   different, the same general principle applies: from a topological homotopy one tries to extract a geometrically well-behaved family of harmonic maps and exploit the geometry of the target space.

More precisely, our main analytic inputs are the harmonic mapping theorem of Schoen--Yau \cite{SchoenYau79} and Noguchi's theory of the space of holomorphic maps to locally symmetric varieties \cite{Nog88}. If two distinct homotopic holomorphic maps exist, these results connect them by a positive-dimensional component $Z$ of the space of holomorphic maps. Then the evaluation map gives a holomorphic totally geodesic immersion $X\times Z\to Y$ after passing to finite covers, and an embedding $\Omega\times\widetilde Z\hookrightarrow\mathcal D$ on universal covers. This converts $h$-deformation rigidity into a geometric obstruction problem. By using conditions on rank, classification results on non-rigid families or Zariski density of monodromy groups, we are able to rule out such a product embedding and establish the $h$-deformation rigidity criteria.

We can already establish certain finiteness result for the equal-rank case.
Combining it with the boundedness part of the Parshin--Arakelov strategy, we ca show the following finiteness theorem.
\begin{thm}(=Theorem \ref{thm:equal-rank-finiteness})
Let $X$ and $Y$ be locally symmetric varieties  , with irreducible source lattice  . If $\operatorname{rk}(X)=\operatorname{rk}(Y)\geq 2$, then $\#\mathrm{Hol}_{\mathrm{nc}}(X,Y)<\infty$.
\end{thm}
We believe the above theorem is known to experts, though we can't find it in the literature. The proof uses the standard ``boundedness plus rigidity'' argument: the boundedness gives a finite-type parameter space after restriction to a curve, while the deformation rigidity shows that each connected component contains at most one map.

At this point an important distinction appears between finiteness and an effective upper bound.  The traditional Parshin--Arakelov strategy works well for proving boundedness and finiteness, but it is usually not well suited to producing \emph{explicit upper bounds}.  In our setting, the boundedness argument restricts maps to a sufficiently general curve in the source space, and places the restrictions in a    finite-type Hom space of maps to the Baily--Borel compactification.  This proves finiteness once rigidity is known. However, this gives little control on the degree of the auxiliary curve or on the number of connected components of the Hom scheme.  Thus a different input is needed if one wants an actual numerical bound.

For an explicit upper bound, we use a different approach.  The basic example is the Kummer map $\tau_K:\mathcal A_2\to\mathcal F_4$.   The $h$-deformation rigidity criteria reduces the problem of counting holomorphic maps to counting integral representations of arithmetic groups. In the Kummer case, an explicit group-ring computation plus some lattice theoretic arguments gives  an   upper bound of the cardinality of such integral representations.  We obtain

\begin{thm}(=Theorem \ref{counterexam_kummer}, Theorem \ref{thm:explicit-kummer-bound})\label{mainthm}
\begin{enumerate}
   \item There exists at least one non-constant holomorphic map from $\mathcal{A}_2$ to $\mathcal{F}_4$ distinct from the Kummer map, modulo $\aut_{orb}(\mathcal{A}_2)$ and $\aut_{orb}(\mathcal{F}_4)$ actions.

    \item $  \#\mathrm{Hol}_{\mathrm{nc}}(\mathcal A_2,\mathcal F_4)<2^{20249}.$
\end{enumerate}

\end{thm}

 Theorem \ref{mainthm} shows that $h$-deformation rigidity should not be confused with global rigidity.  
 This gives a negative answer to the strongest form of the rigidity expectation suggested in \cite{Farb24}. This upper bound is certainly very far from optimal.  But  it is completely \emph{explicit}, which illustrates the role of $h$-deformation rigidity in effective Shafarevich-type problems: it allows one to pass from holomorphic maps to discrete arithmetic data which can, at least in principle, be counted.

We also obtain several supplementary results which may be useful for further examples. In Sect.~\ref{alternative_finite}, we give an alternative proof of the finiteness of $\mathrm{Hol}_{\mathrm{nc}}(\mathcal A_2,\mathcal F_4)$ using Hermitian metric rigidity, the classification of totally geodesic embeddings, and the Borel--Harish-Chandra theorem (see the proof of Theorem~\ref{thm:alternative-finiteness}).  We analyze orbifold automorphism groups of the orthogonal quotients appearing in the Kummer construction, and in Sect.~\ref{sec4} we compute the orbifold automorphism group of the Eisenstein--Picard modular $n$-fold considered there  , for even $n\geq2$,   as the cyclic group $C_{\gcd(6,n+1)}$ (see Proposition~\ref{Eisenstein–Picard}).  We also discuss period maps attached to low-genus curves, cubic hypersurfaces, configurations of points, and Albanese varieties.

A particularly interesting next case is the intermediate-Jacobian period map for smooth cubic threefolds
\[
  \mathcal C_3\longrightarrow \mathcal A_5.
\]
Its $h$-deformation rigidity is already covered by Theorem~\ref{h-deform-rigidity_thm_intJac_case}. By Torelli theorem, $\mathcal C_3\subset \bB^{10}/\Gamma$ and the complement consisting of two    Heegner   divisors. Following \cite{BFJS26} we show that $\aut_{orb}(\bB^{10}/\Gamma)=\{1 \}$.  Thus the remaining problems are of the type highlighted above: determine whether $\aut_{orb}(\cC_3)=\{1\}$,  understand the possible nontrivial group homomorphisms    $\pi_1^{orb}(\cC_3)\to \Sp(10,\mathbb Z)$  , and if global rigidity fails, seek an explicit upper bound for the number of holomorphic maps.  We plan to investigate this case in future work.

\subsection*{Organization of the paper}
In \S\ref{sec2}, we collect basic results on bounded symmetric domains and lattice theory and recall the period maps studied in the paper.  In \S\ref{sec3}, we introduce $h$-deformation rigidity for orbifold maps, prove the extension theorem using harmonic mapping theory, and establish the rigidity criteria and applications above.  In \S\ref{sec:equal-rank-finiteness}, we combine $h$-deformation rigidity with the Parshin--Arakelov boundedness method to prove finiteness in the equal-rank case.  We then show the failure of global rigidity for the Kummer construction and give an alternative arithmetic proof of finiteness.  In \S\ref{sec:explicit-bound}, we carry out the effective integral and lattice-theoretic counting which gives the bound $2^{20249}$.  Finally, in \S\ref{sec4}, we discuss further examples and open problems, including the cubic-threefold case.

\subsection*{Acknowledgements} We would like to thank Runlin Zhang for helpful conversations.  C. D. is partially  supported by the National Natural Science Foundation of China
(Grant No. 12601146), start-up funding from Shenzhen University and the Shenzhen Peacock Plan. F.S. is partially supported by the Fundamental Research Funds for the Central Universities,  Shaanxi Fundamental Science Research Project for Mathematics and Physics (Grant No. 25JSZ007,25JSQ045) and Shaanxi NSF (No. 2025JC-QYCX-002). R.S. is supported by the National Natural Science Foundation of China (Grant No. 12501059) and a grant from Xiamen University (No. X2450215).

\subsection*{AI disclosure} Our project starts from September 2025.  We obtain all the main results except for an explicit upper bound of  $ \mathrm{Hol}_{\mathrm{nc}}(\cA_2,\cF_4)$ before February 2026. Then 
Section \ref{sec:explicit-bound} of this article was obtained through discussions with ChatGPT Astra. We also use Claude Fable in revising the arguments. The authors take full responsibility for the    mathematical content and references  .

\section{Preliminary}\label{sec2}
\subsection{Basics on quotient of bounded symmetric domains}

\begin{defn}
    A bounded symmetric domain $\cD$ is a  non-empty bounded open subset in $\bC^n$ such that  the  group $\mathrm{Hol}(\cD)$ of holomorphic automorphisms acts  transitively and there is a point    $p\in\cD$ which is an isolated fixed point of a holomorphic involution $s_p:\cD\to\cD$  . 
\end{defn}
By the classification,  there are four classical families and two exceptional families.
\begin{exmp}(Classical domains, cf. {\cite[Chapter 4]{Mok89}})
Let $M(p,q)$ be the space of $p \times q$ complex    matrices   and $A^\ast$ denote the    conjugate transpose of   $A$. Then
\begin{equation*}
    D^I_{p,q}:=\{ A \in M(p,q)\ |\ \mathrm{Id}-A^\ast A>0\}, p\ge q\ge 1
\end{equation*}
  is the type I bounded symmetric domain.  The type II, III and IV bounded symmetric domains are defined as follows
\begin{equation*}
     \begin{split}
        &D^{II}_n:=\{ A\in D^I_{n,n} \ |\ A^t=-A  \} ,\ \ \ n\ge 2; \\
        &D^{III}_n:=\{ A\in D^I_{n,n} \ |\ A^t=A  \} ,\ \ \ n\ge 1; \\
        &D^{IV}_n:=\{ z\in \bC^n\ |\ |z|^2<2,\ |z|^2<|\frac{\mathop{\sum}z_i^2}{2}|^2+1 \}.
  \end{split}  
\end{equation*}
\end{exmp}

As a Hermitian symmetric space of non-compact type, the bounded symmetric domain is a homogeneous space $G/K$ where $G$  is a non-compact real semisimple Lie group with trivial center and $K\le G$ is its maximal compact subgroup (cf. {\cite[Chapter 3 \S 1]{Mok89}}). Moreover, $G/K$ has seminegative curvature, cf. \cite[Proposition 2, p46]{Mok89}.

\begin{defn}
    The rank of a bounded symmetric domain $\cD$ is the maximum dimension of a flat and totally geodesic submanifold of $\cD$.  The rank of a locally symmetric variety is defined to be the rank of its universal covering domain.
\end{defn}
\begin{rem}
Let $\mathfrak{g}$ (resp. $\mathfrak{k}$ ) be the Lie algebras of $G$ (resp. K) and we have the Cartan decomposition $\mathfrak{g}=\mathfrak{k}\oplus \mathfrak{m}$.
An equivalent definition of rank for $\mathcal{D}=G/K$ is defined as the dimension of a maximal abelian subspace of $\fm$, see \cite[Appendix \S 2]{Mok89}. In particular, $\rank(\cD) \le \dim \mathfrak{h}$    where $\mathfrak{h}\subset \mathfrak{g}$ is the Cartan subalgebra.
\end{rem}
The ranks of four classical domains are given by the following table. 
\begin{table}[htbp]
\centering
\footnotesize
\caption{rank of classical domains}
\begin{tabular}{|c|c|c|c|c|} 
\hline  
 & $D^I_{p,q}$  & $D^{II}_n$ & $D^{III}_n$ & $D^{IV}_n$ \\
\hline
$G/K$ & $\mathrm{SU}(p,q)/\mathrm{S}(\mathrm{U}(p) \times \mathrm{U}(q)) $ & $\SO^\ast(2n)/U(n)$ & $\Sp_{2n}(\bR)/U(n)$  &    $\SO_0(n,2)/(\mathrm{SO}(n)\times\mathrm{SO}(2))$   \\
\hline
rank & $q$ & $[\frac{n}{2}]$ &$n$ & $2$ \\
\hline
\end{tabular}
\label{tab:summary}
\end{table}

\begin{thm}(The Polydisc Theorem, cf. {\cite[Chapter 5 \S 1]{Mok89}}) \label{Theorem:PD}
  Let $(G/K,g)$ be a    Hermitian symmetric space of noncompact type and   rank $r$,   with the invariant metric normalized on its maximal polydiscs,   then there exists a totally geodesic complex    submanifold   $\Delta$ such that $(\Delta,\Delta|_g)$ is isometric to a Poincar\'e polydisc\[ \Big(\{z=(z_1,\cdots,z_r)\in \bC^r\ |\ |z_i| <1 \},\ ds^2=\mathop{\sum} \limits_{i=1}^r \frac{4|dz_i|^2}{(1-|z_i|^2)^2} \Big).\]
\end{thm}

Let $\Gamma\leq\aut(\mathcal D)$ be a torsion free lattice and  $X_\Gamma=\mathcal D/\Gamma$. We recall the following theorem as a corollary of the Hermitian metric rigidity theorem due to Mok and To \cite{Mok89,To89}. This result also applies to orbifold maps.
\begin{thm}\label{thm:MR}
Assume $\operatorname{rk}(\mathcal D)\geq2$. Let $f:X_\Gamma\to N$ be a nonconstant holomorphic orbifold map.
If $\mathcal D$ is irreducible and $(N,h)$ is K\"{a}hler with nonpositive curvature in the sense of Griffiths, then $f$ is an orbifold immersion and is totally geodesic.
\end{thm}

\subsection{Basics on lattice theory}
We recall some basic notions and facts in lattice theory. Our main references are \cite[Chapter 14]{Huybrechts2016} and \cite{Niku79}. A lattice $L$ is a free $\mathbb{Z}$-module of finite rank with a symmetric non-degenerate bilinear form $(,)$. \textit{Signature} of a lattice is defined as the signature of the $\mathbb{R}$-extension of $(,)$. $L$ is called an even lattice if $(x,x)\in 2\mathbb{Z}$ for any $x\in L$. The determinant of the intersection matrix with respect to arbitrary basis is called the \textit{discriminant}, denoted by $\operatorname{disc}(T)$. Let $L^\vee:=\operatorname{Hom}_\mathbb{Z}(L, \mathbb{Z})$ be the dual lattice. We have natural inclusions $L\hookrightarrow L^\vee \hookrightarrow L_\mathbb{Q}$. The cokernel of the inclusion $L\hookrightarrow L^\vee$ is called the discriminant group $A_L:=L^\vee/L$, which is a finite group of order    $|A_L|=|\operatorname{disc}(L)|$  . $L$ is called \textit{unimodular} if $A_L$ is trivial. For a given lattice $L$, its twist $L(m)$ is given by changing the intersection form as $(,)_{L(m)}:=m(,)$.
\begin{exmp}
    The hyperbolic lattice $U$ is the lattice given by the intersection matrix \[\begin{pmatrix}
        0 & 1\\
        1 & 0
    \end{pmatrix}.\]
    The $E_8$-lattice, denoted by $E_8$, is the lattice with intersection matrix the same as the $E_8$-Cartan matrix.
\end{exmp}

   For an even lattice   $L$   , the discriminant quadratic form is $q_L:A_L\to\mathbb Q/2\mathbb Z$, $q_L(\overline x)=(x,x)\bmod2\mathbb Z$; its associated bilinear form takes values in $\mathbb Q/\mathbb Z$  . The finite group $A_L$ together with $q_L$ is called the \textit{discriminant form} of $L$. A morphism between two lattices $L_1,L$ is a linear map that respects the
quadratic forms. An injective morphism between two lattices $L_1\hookrightarrow L$ is called a \textit{primitive embedding} if $L/L_1$ is torsion free, equivalently, $L_1=L\cap (L_1)_\mathbb{Q}$.  We collect two basic propositions here, which must be familiar to experts.

\begin{pro}\label{finite_index_disc}
    If $L_1 \hookrightarrow L$ has finite index, then $\operatorname{disc}(L_1)=\operatorname{disc}(L)(L:L_1)^2$.
\end{pro}

\begin{pro}(\cite[Proposition 1.5.1]{Niku79})\label{orthogonal_disc_iso}
    A primitive embedding $L_1\hookrightarrow L$ into an even lattice with $L_2:=L_1^\perp$ is determined by subgroups    $C_{L_1}\subset A_{L_1}, C_{L_2}\subset A_{L_2}$   with an isomorphism $(C_{L_1}, q_{L_1}|_{C_1})\cong (C_{L_2}, -q_{L_2}|_{C_2})$. In particular, if $L$ is unimodular, the primitive embedding is determined by an isomorphism $(A_{L_1}, q_{L_1})\cong (A_{L_2}, -q_{L_2})$. 
\end{pro}
\begin{proof}
    For a fixed even lattice $L'$, there is a bijective correspondence between finite index overlattice $L'\subset L$ and isotropic subgroups $H\subset A_{L'}$. Namely, for a finite index overlattice $L'\subset L$, we take $H=L/L' \subset L^\vee/L'\subset (L')^\vee/L'=A_{L'}$. Conversely, for a given subgroup $H \subset A_{L'}$, we take $L$ to be its pre-image under $(L')^\vee \rightarrow A_{L'}$.

  For a primitive embedding $L_1 \hookrightarrow L$ into an even lattice with $L_2:=L^\perp_1$, one associates the  finite index sublattice $L':=L_1\oplus L_2 \subset L$. The associated isotropic subgroup is    $H=L/(L_1\oplus L_2) \subset A_{L'}\cong A_{L_1}\oplus A_{L_2}$  . The two projections
 \[H\longrightarrow A_{L_i},\]
    $i=1,2$, are injective, as both inclusions $L_1,L_2\subset L$ are primitive. Let $C_1, C_2$ denote the images of the projections respectively. Then $C_1 \cong H \cong C_2$ gives the anti-isometry.

 In contrast, given subgroups    $C_{L_1}\subset A_{L_1}, C_{L_2}\subset A_{L_2}$   with isometry $(C_{L_1}, q_{L_1}|_{C_1})\cong (C_{L_2}, -q_{L_2}|_{C_2})$, the isometry graph gives an isotropic subgroup $H$, which gives rise to a finite index overlattice $L_1\oplus L_2 \subset L$. As $H$ is a graph of an isomorphism, the embedding $L_1 \hookrightarrow L$ must be primitive. 

 If $L$ is unimodular, from Proposition \ref{finite_index_disc}, we know $|H|^2=|A_{L_1}||A_{L_2}|$, which forces $C_{L_1}=A_{L_1}, C_{L_2}=A_{L_2}$. Hence the conclusion follows. 
\end{proof}

\subsection{Construction of period maps}
 \subsubsection{Kummer construction}\label{Kummer}
Let $A$ be a principally polarized abelian surface and $\iota: A \xrightarrow{} A$ the natural involution. Then $A/\iota$ has $16$ isolated  singularities locally isomorphic to $\{x^2+y^2+z^2=0\}$. Then there is a crepant resolution \[ S \xrightarrow{} A/\iota \]
such that  $S$ is a K3 surface, which is known as Kummer K3.   The Neron-Severi group  of a Kummer K3 is an even  lattice with rank $\ge 17$, 
e.g., see \cite{Huybrechts2016} for more details.  Given  a sublattice $N\subset\NS(S)$ with signature $(1,r)$ such that $N \hookrightarrow\NS(S)$ is primitive and $N$ contains  a nef and big class, denote $N^\perp$ the orthogonal complement of $$N \hookrightarrow L_{K3}=H^2(S,\bZ)\cong U^3\bigoplus E_8^2(-1).$$ Then $T:=N^\perp$ is an even lattice of signature $(2,19-r)$. 

Let $ \fX \to \bD^\circ$ be a family of principally polarized abelian surfaces, which gives a monodromy operator $\gamma\in \Sp(4)$. Then the above Kummer K3's construction can be done relatively over $\bD^\circ$ and it turns out that there is a family of Kummer K3 surfaces $\pi: \cY \to \bD^\circ$ associated to $ \fX \to \bD^\circ$. Then $R\pi_\ast^2 \bZ $ is a polarizable Hodge structure of weight $2$. Once we fix an embedding $N \hookrightarrow\NS(\cY_t)$, then parallel transport will ensure that the monodromy representation $\rho: \pi_1(\bD^\circ,t) \to O(H^2(\cY_t,\bZ))$ satisfy $\rho(1)(N)=N$, e.g.,see ~\cite[\S 4]{DHNT15}. Thus we have a monodromy $\rho_N: \Sp(4) \to O(L_{K3})$ such that $\rho_N(\gamma)(N)=N$. Then the  construction produces a holomorphic map
\[  p_N:\ \cA_2 \rightarrow \cD_{N}/\Gamma_{N} \]
known as period map, where $$\cD_{N}:=\{ z\in \bP(T \otimes \bC)\ |\ z.z=0,\ z\overline{z}>0\}$$ is the period domain (isomorphic to the Type IV classical domain $D_{19-r}^{IV}$) , and $\Gamma_{N} \le O(T)$ 
is the image of 
\[ \{ \gamma\in O(L_{K3})\ |\ \gamma(N)=N\} \to O(T),\ \gamma\mapsto \gamma|_{T}.\]

Moreover, if $N$ is monodromy invariant, i.e,. the monodromy operator $\gamma$ fixes $N$ pointwise,  then by the lattice gluing trick, 
 \[\Gamma_{N}=\widetilde{SO}(T):=\{\gamma \in SO(T)\ |\ \gamma=\mathrm{Id}\  \hbox{on} \ \ T^\vee/T\}\] is the discriminant kernel in $\SO(T)$.

It is known that  there is a monodromy invariant, nef and big class $h$ with $h^2=4$, for example see \cite{Huybrechts2016}. Then taking $N_0=\langle h \rangle$,   we have the period map \[p_0: \mathcal{A}_2 \rightarrow \mathcal{F}_4:=\mathcal{D}_{N_0}/\Gamma_{N_0}\] 
to the moduli space of degree $4$ quasi-polarized K3 surfaces.
Note that if there are embeddings $N_0\subset N_1 $ where are all monodromy invariant, it will induce the natural embedding   $ \cD_{N_1} \hookrightarrow \cD_{N_0} $  and       also \[ \{ \gamma\in O(L_{K3})\ |\ \gamma_{N_0}=\mathrm{Id}\} \to \{ \gamma\in O(L_{K3})\ |\ \gamma_{N_1}=\mathrm{Id}\}. \]
Then in this case the period maps $p_{N_0}:\ \mathcal{A}_2  \xrightarrow{} \cD_{N_0}/\Gamma_{N_0}$ will factor through 
\begin{equation}
    \mathcal{A}_2  \xrightarrow{p_{N_1}} \cD_{N_1}/\Gamma_{N_1}  \xrightarrow{\iota}  \cD_{N_0}/\Gamma_{N_0}
\end{equation} where $\iota: \cD_{N_1}/\Gamma_{N_1}  \xrightarrow{}  \cD_{N_0}/\Gamma_{N_0}$ is induced by $ \cD_{N_1} \hookrightarrow \cD_{N_0} $ .

\begin{rem}
  Here we start with a principally polarized abelian surface.   If we consider a general   polarized abelian surface $A$ of polarization type $(1,d)$, there is also a period map $p_{d,T}:\ \cA_{2,d} \rightarrow \cD_T/\Gamma $ where $\cA_{2,d}$ is the coarse moduli space of the  abelian surfaces  of polarization type $(1,d)$. In this case,  $\cA_{2,d}\cong \bH_2/\Gamma_d$  where (cf. \cite[\S 1]{Sankaran97})$$\Gamma_d=\Big\{ \gamma\in \Sp(4,\bQ) \ |\ \gamma\in \begin{pmatrix}
\bZ  & \bZ & \bZ & d\bZ \\
d\bZ  & \bZ & d\bZ & d\bZ \\
\bZ  & \bZ & \bZ & d\bZ \\
\bZ  & \frac{1}{d}\bZ & \bZ & \bZ 
\end{pmatrix} \Big \}.$$  
\end{rem}

\subsubsection{Hypersurfaces with interesting intermediate Jacobian} \label{period:hypersurface}
Suppose $d\ge 3$ and $n\ge 2$ or $(n,d)=(1,d)$ for $d\ge 4$.
Let \[U_{n,d}\subset |\cO_{\bP^{n+1}}(d)|\] be 
the locus of smooth degree $d$ hypersurfaces in $\bP^{n+1}$. Denote $M_{n,d}$  the good moduli space of the quotient stack $[U
_{n,d}/\PGL(n+2)]$. It is well known that a smooth hypersurface of dimension $n\ge 2$ and degree $d\ge 3$ in complex projective space  has finite projective.    Thus $M_{n,d}$ is an  orbifold.
By \cite{Collino86}, it is known that only for \[ (n,d) \in \{ (1,d), (3,3),(3,4), (5,3)   \}, \]
the intermediate Jacobian
\[  J_X= \frac{H^n(X,\bC)}{F^{\frac{n+1}{2}} H^n(X,\bC)\oplus H^n(X,\bZ) }\]
is a p.p.a.v. Then in these cases, there is a period map 
\[ p:\ M_{n,d} \rightarrow \cA_{N(n,d)},\ [X] \ \mapsto\ [J_X]\]
where $N(n,d)$ is the dimension of the abelian variety $J_X$. See the following table for precise values of $N(n,d)$.

\begin{table}[htbp]
\centering
\caption{dimension of intermediate Jacobian}
\begin{tabular}{|c|c|c|c|c|} 
\hline  
$(n,d)$ &  $(1,d)$ & $(3,3)$  & $(3,4)$  &  $(5,3)$    \\
\hline
$N(n,d)$ &  $\frac{(d-1)(d-2)}{2}$ & $5$ & $30$  &  $21$ \\
\hline
\end{tabular}
\label{tab:dim}
\end{table}

These moduli spaces are well studied in algebraic geometry. For example, by the cyclic covering $Y \xrightarrow{} \bP^4$ branched along $X$ and the global Torelli theorem, there is an open immersion \[ M_{3,3} \hookrightarrow \bB^{10}/\Gamma  \]
where $\Gamma$ is the image of arithmetic subgroup $\mathrm{SU}(\cE,h)$ under the projection $\mathrm{SU}(10,1) \to \mathrm{PU}(10,1)$. The complement $\bB^{10}/\Gamma-M_{3,3}$ consists of two    Heegner   divisors, see \cite[Theorem 6.1]{ACT11} for more details.

Note that for $n=1$, the moduli space $M_{1,d}$ is the moduli space of plane curves. In this case we denote $g(d)=N(1,d)=\frac{(d-1)(d-2)}{2}\ge 3$ (as $d\ge 4$) and $ \cM_{g}$ the moduli    space   of smooth curves of genus $g$. Then there is factorization \[  p=\tau \circ \theta  \]  for period map $p:M_{1,d} \rightarrow \cA_{g(d)}$ where $\theta:\ M_{1,d} \rightarrow \cM_{g(d)}$ is the moduli map and $\tau: \cM_{g(d)} \rightarrow   \cA_{g(d)}$  is the Torelli morphism.
\begin{rem}
    The above construction of period for moduli space of plane curves can be generalized to the moduli  space of smooth curves $U\subset |L|$ in a linear system of line bundle $L$ on a smooth projective surface $S$. 
\end{rem}

\subsubsection{ Intersections of two quadrics of odd dimension} \label{period:ci}
Suppose $X$ is a smooth complete intersection of two quadrics in $\bP^{2n+3}$ for $n\ge 0$. Then Hodge decomposition of the middle cohomology  is of the form $H^{2n+1}(X,\bC)=H^{n,n+1}(X) \oplus H^{n+1,n}(X)$ . Then its intermediate Jacobian \[ J_X:=\frac{H^{n,n+1}(X,\bC)}{H_{2n+1}(X,\bZ)}\]
is a p.p.a.v of dimension $(n+1)$, see \cite{Reid72}.  Let $U \subset \Gr(2,H^0(\cO_{\bP^{2n+3}}(2)))$ be the locus of smooth complete intersections and $M_n$ the good moduli space of the quotient stack $[U/\PGL(2n+4)]$, then there is a period map 
\begin{equation*}
    p:\ M_n \rightarrow \cA_{n+1}.
\end{equation*}

\subsubsection{ Moduli space of  points on affine line }\label{Periodmap:points}
Let $\mathrm{UC}_n=(\bC ^n-\Delta_n)/S_n$ be the space of unordered $n$-points on $\bC$ where $\Delta_n:=\{ (\lambda_1,\cdots,\lambda_n)\in \bC^n |\ \lambda_i= \lambda_j \ \hbox{for some}\ i\neq j\  \}$ and $S_n$ acts on $\bC ^n-\Delta_n$ by standard permutation of $n$ points. Denote $\cD_{2l}/\Gamma_{2l}$ the locally symmetric variety associated to quasi-polarized K3 surface of degree $2l$. By \cite[\S 12,Examples]{DK06}, we have the following period mapps.
\begin{enumerate}
    \item $n=5$, the distinct $5$ points $(\lambda_1,\cdots,\lambda_5)\in \mathrm{UC}_n$ can naturally define a plane  sextic curve $C$ by \[  z_0^6=z_0(z_1-\lambda_1z_2)\cdots (z_1-\lambda_5z_2).\]
    The double cover of $\bP^2$ branched along $C$ is a K3 surface with $A_n$-type singularities at worst, whose minimal resolution is a smooth K3 with a  quasi-polarization  of degree $2$. Thus there is a period map 
    \begin{equation}
        \mathrm{UC}   _5  \rightarrow\cD_2/   \Gamma  _2
     \end{equation}
    \item $n=6$, given  distinct $6$ points $(\lambda_1,\cdots,\lambda_6)\in \mathrm{UC}_n$, there is a natural degree $6$ hypersurface \[ X=\{ z_3^3+z_2^3+(z_1-\lambda_1z_0)\cdots (z_1-\lambda_6z_0)=0 \} \subset \bP(1,1,2,2)\] in weighted projective $3$-fold $\bP(1,1,2,2)$. Then $X$ is a K3 surface with three $A_1$-singularities.  Then its minimal resolution is a smooth K3 with a    quasi-polarization   of degree $6$. This construction also produce a period map 
    \begin{equation}
        \mathrm{UC}_n \rightarrow \cD_6/\Gamma_6
    \end{equation} as before. 
    \item $n=8$,    use the finite cover $\mathrm{UC}_{4,4}\to\mathrm{UC}_8$ marking a partition into two sets of four points. Let   $([z_0,z_1],[w_0,w_1])$ be bi-homogeneous coordinates on $\bP^1 \times \bP^1$, then any distinct $8$ points $(\lambda_1,\cdots,\lambda_8)\in \mathrm{UC}_n$ produce a bidgree $(4,2)$-divisor $C$ on $\bP^1 \times \bP^1$ by \[ w_0^2(z_1-\lambda_1z_0)\cdots (z_1-\lambda_4z_0)+w_1^2(z_1-\lambda_5z_0)\cdots (z_1-\lambda_8z_0)=0.  \]
    Then the double cover of $\bP^1 \times \bP^1$ branched along the divisor $C+\{w_0=0\}+\{w_1=0\}$ is a K3 with eight  $A_1$-singularities. Then there is a period map 
    \begin{equation}
        \mathrm{UC}   _{4,4}  \rightarrow\cD_4/\Gamma_4
    \end{equation} as constructed before. 
\end{enumerate}

\section{Deformation rigidity of period maps} \label{sec3}
\begin{defn}[$h$-deformation rigidity]
Let $X$ and $Y$ be two complex analytic spaces. Let $f:\,X \to Y$ be a holomorphic map. We say $f$ is {\rm $h$-deformation rigid} if for every holomorphic map $g:\,X \to Y$ such that $f$ is homotopic to $g$, then $g=f$.
\end{defn}
\begin{rem}
   $h$-deformation rigidity implies the deformation rigidity. 
\end{rem}
In the context of locally symmetric varieties and moduli spaces, one usually needs to consider the orbifold structure. Let $X$ and $Y$ be two complex orbifolds (i.e. properly discontinuous quotients of complex manifolds). Consider a holomorphic map $f:\,X \to Y$ respecting the orbifold structures. Then it makes sense to discuss the homotopy class of continuous (or holomorphic) maps in the category of orbifolds (cf. {\cite[Remark~2.1]{Farb22}}). Therefore, we also have the notion of $h$-deformation rigidity for holomorphic maps between complex orbifolds.

\subsection{An extension theorem}

Let $X:=\Omega/\Gamma_1$, $Y:=\mathcal D/\Gamma_2$ be two locally symmetric varieties, namely arithmetic quotients of bounded symmetric domains. Note that $X$ and $Y$ are naturally complex orbifolds.
\begin{thm}[Schoen-Yau]\label{SchYau}
  Let $X$ be a smooth quasi-projective variety equipped with a complete K\"ahler metric $h_X$ with finite volume. Let $Y$ be a locally symmetric variety with nonpositive sectional curvatures, and let $h_Y$ be the invariant metric on $Y$. Let $f,g:\,X \to Y$ be homotopic harmonic maps of finite energy. Then there exists a harmonic map
  \[
\Psi:\, X \times [0,1] \to Y
  \]
  with respect to the product metric    $h_X \oplus dt^2$   on $X \times [0,1]$ ($dt^2$ is the Euclidean metric on $[0,1]$), such that
  \begin{itemize}
    \item[i)] $f=\Psi|_{X \times \{0\}}$, $g=\Psi|_{X \times \{1\}}$ and $\Psi$ provides a homotopy between $f$ and $g$, equivalent to the given one.
    \item[ii)] For each $x \in X$, the curve $\gamma_x:\,[0,1] \to Y$, $t \mapsto \Psi(x,t)$ is a geodesic in $Y$ with respect to the invariant metric $h_Y$, and the energy of $\Psi(\cdot, t)$ is constant in $t$.
  \end{itemize}
\end{thm}

Next we will use the moduli space of holomorphic maps. Let $N$ be a
smooth connected quasi-projective variety and let
$Y=\mathcal D/\Gamma$ be a smooth locally symmetric variety.
Denote by $\mathrm{Hol}(N,Y)$ the space of holomorphic maps from $N$ to $Y$,
endowed with the compact-open topology. The following statement
combines \cite[Theorem~4.7]{Nog88} with
\cite[Theorem~2.12 and Lemma~2.13]{MiyanoNoguchi91}.

\begin{thm}[Noguchi]\label{Nog-thm}
\begin{itemize}
 \item[i)] $\mathrm{Hol}(N,Y)$ is a smooth complex analytic space,
 and each connected component is quasi-projective.
 \item[ii)] For each connected component
 $Z\subset\mathrm{Hol}(N,Y)$ and a point $x\in N$, the evaluation map
 \[
 e_x:Z\longrightarrow Y,\qquad u\longmapsto u(x)
 \]
 is a proper holomorphic immersion onto a totally geodesic submanifold.
 The metric $h_Z=e_x^*h_Y$ is complete and independent of $x$.
\end{itemize}
\end{thm}

\begin{thm}\label{Ext-thm}
Suppose that $\Gamma_1$ is irreducible and the rank of $\Omega$ is greater
than one. Let $f:X\to Y$ be a nonconstant holomorphic map of orbifolds.
If there exists another holomorphic map $g:X\to Y$ of orbifolds, homotopic
to $f$ but distinct from it, then $f$ and $g$ can be enlarged in the
following way: after replacing $X,Y$ by compatible finite orbifold
covers, there exist a positive-dimensional connected smooth
quasi-projective variety $Z$, with a complete locally symmetric K\"ahler
metric, and a holomorphic totally geodesic immersion
\[
 \Phi:X\times Z\longrightarrow Y
\]
for a product metric, such that $f$ and $g$ are its restrictions
to $X\times\{z_1\}$ and $X\times\{z_2\}$ for two points $z_1,z_2\in Z$.
On universal covers, the induced map
$\widetilde\Phi:\Omega\times\widetilde Z\hookrightarrow\mathcal D$
is a totally geodesic embedding.
\end{thm}

\begin{proof}
By definition, the homotopy between $f$ and $g$ respects the orbifold
structures. One can replace $\Gamma_i$ by finite-index subgroups and
assume that $X,Y$ are torsion-free quotients and the homotopy
lifts. By \cite{BB66}, $X,Y$ are quasi-projective manifolds. The maps
$f,g$ are harmonic and have finite energy for the complete
finite-volume K\"ahler metric of \cite[Lemma~4.2]{Nog88}.
By Theorem~\ref{SchYau}, we can connect $f$ and $g$ by a one-parameter
family of harmonic maps. From \cite[Lemma~4.6]{Nog88}, each
$\Psi_t:X\to Y$, $0\leq t\leq1$, is holomorphic.
There is the evaluation map
\[
 \Phi:X\times\mathrm{Hol}(X,Y)\longrightarrow Y,\qquad (x,u)\longmapsto u(x).
\]
The holomorphic maps $\Psi_t$ define a continuous path
$\ell:[0,1]\to\mathrm{Hol}(X,Y)$, and $\Psi$ is the composite
\[
 X\times[0,1]
 \xrightarrow{(\mathrm{Id},\ell)}X\times\mathrm{Hol}(X,Y)
 \xrightarrow{\Phi}Y.
\]
Let $Z\subset\mathrm{Hol}(X,Y)$ be the connected component containing
the image of $\ell$. In particular, $f=\ell(0)$ and $g=\ell(1)$ belong
to $Z$. Since $f\ne g$, there is a point $x\in X$ with $f(x)\ne g(x)$.
The geodesic $\gamma_x$ connects these two points and is contained in
the image of $Z$ under $\Phi(x,\cdot)$. Thus $\dim Z>0$, and
we have the following diagram:
\[
 \xymatrix{
 \{x\}\times[0,1]\ar@{^{(}->}[d]\ar[rd]^{\gamma_x}&\\
 X\times[0,1]\ar[d]_{(\mathrm{Id},\ell)}\ar[r]^{\Psi}&Y\\
 X\times Z\ar@{^{(}->}[r]&X\times\mathrm{Hol}(X,Y).\ar[u]_{\Phi}
 }
\]
Since $f$ is nonconstant and the energy of $\Psi_t$ is constant in $t$,
each $\Psi_t$ is nonconstant. By Theorem~\ref{thm:MR}, each is a totally
geodesic immersion. Every element of $Z$ is homotopic
to $f$, so the same argument applies to all the maps parametrized by $Z$.

To obtain total geodesity of the product, we consider the maps
$e_x=\Phi(x,\cdot):Z\to Y$. They belong to one connected component
$W\subset\mathrm{Hol}(Z,Y)$ and define a holomorphic map
\[
 j:X\longrightarrow W,\qquad x\longmapsto e_x.
\]
The theorem about evaluation map of Miyano--Noguchi
\cite[Theorems~2.15 and 4.9]{MiyanoNoguchi91} gives a totally geodesic
isometric immersion
\[
 E:Z\times W\longrightarrow Y,\qquad (u,w)\longmapsto w(u),
\]
where $Z\times W$ carries the product of the complete evaluation metrics.
Note that the map $e_f=E(f,\cdot):W\to Y$ is totally geodesic and
$e_f\circ j=f$. Put $h_X=f^*h_Y$. Then by Theorem~\ref{thm:MR}, this metric
lifts to $c \cdot g$ on $\Omega$, where
$g$ is the invariant metric on $\Omega$ and $c>0$ a constant.
Since $de_f$ is injective, the identity
\[
 0=\nabla df=de_f\circ\nabla dj
\]
shows that $j$ is totally geodesic. Moreover, $j^*h_W=h_X$.
Consequently $\Phi(x,u)=E(u,j(x))$ is a totally geodesic immersion with
\[
 \Phi^*h_Y=h_X\oplus h_Z.
\]
Its restrictions at $f,g\in Z$ are the two prescribed maps.

The lift
$\widetilde\Phi:\Omega\times\widetilde Z\to\mathcal D$ is therefore
a totally geodesic isometric immersion between complete simply
connected manifolds of nonpositive sectional curvature.
\end{proof}

\subsection{$h$-deformation rigidity criteria}
Now we will use Theorem~\ref{Ext-thm} and Theorem~\ref{SchYau} in the
previous subsection to derive $h$-deformation rigidity of certain
holomorphic maps.

\begin{thm}[Equal rank case]\label{h-deform-rigidity_thmER}
Let $X,Y$ be locally symmetric varieties of rank greater than one,
with $X$ irreducible. Suppose that the rank of $X$
equals the rank of $Y$. Then any nonconstant holomorphic map
$f:X\to Y$ of complex orbifolds is $h$-deformation rigid.

In particular, the period map
\[
 p_T:\mathcal A_2\longrightarrow\mathcal D/\Gamma_T
\]
induced by the Kummer construction in Sect.~\ref{Kummer} is
$h$-deformation rigid.
\end{thm}

\begin{proof}
If $f$ fails to be $h$-deformation rigid, Theorem~\ref{Ext-thm} gives
a totally geodesic embedding
$\Omega\times\widetilde Z\hookrightarrow\mathcal D$, with $\dim Z>0$.
A maximal flat in $\Omega$, multiplied by a complete geodesic in
$\widetilde Z$, gives a flat of dimension $\operatorname{rk}(X)+1$
in $\mathcal D$. This contradicts $\operatorname{rk}(X)=\operatorname{rk}(Y)$.
\end{proof}

\begin{thm}[Ball quotients and moduli of K3]\label{h-deform-rigidity_thm_ball_case}
Let $X$ be a finite-volume ball quotient of dimension at least two, and
let $Y=\mathcal D/\Gamma_T$ be as in Sect.~\ref{Kummer}. If a
holomorphic orbifold map $f:X\to Y$ has {\rm rank at least two} at a
general point of $X$, then $f$ is $h$-deformation rigid. 
In particular,
every quasi-finite holomorphic orbifold map $X\to Y$ is
$h$-deformation rigid.
\end{thm}

\begin{proof}
Suppose that $f$ is homotopic to a distinct holomorphic orbifold map
$g:X\to Y$. Pass to compatible smooth finite orbifold covers so that
the homotopy lifts, taking the target group to be a torsion-free normal
subgroup of finite index in $O^+(T)$ contained in $\Gamma_T$. We keep
the notation $X,Y,f,g$.

Note that the ball quotient $X$ is quasi-projective by \cite[\S2.2]{Mok12}. One can choose a complete finite-volume K\"ahler metric on $X$ as in
\cite[\S4]{Nog88}, for which $f$ and $g$ have finite energy
\cite[Lemma~4.2]{Nog88}. We then use Theorem~\ref{SchYau} of Schoen--Yau to join
$f$ to $g$ by a one-parameter family of harmonic maps. By Noguchi's
lemma \cite[Lemma~4.6]{Nog88}, each map in this family is holomorphic.
Thus the connected component $Z$ of $\mathrm{Hol}(X,Y)$ containing $f$
has positive dimension, and the evaluation map gives us a nontrivial holomorphic
deformation
\[
  \Phi:X\times Z\longrightarrow Y,\qquad (x,u)\longmapsto u(x).
\]

Now consider $\mathbb V$, the pullback of the tautological polarized $\mathbb Z$-VHS of K3 type on $Y$ by $f$. It is
non-isotrivial because $f$ is nonconstant. We then obtain a nontrivial deformation of $\mathbb V$ by pulling back along $\Phi$ (cf. \cite[\S3.3]{Saito-Zucker}).

By the classification results of Saito--Zucker
\cite{Saito-Zucker} (see in particular Theorem~5.2.3 and eq.~(5.3.7)--(5.3.10) in \emph{loc.cit.}),
the flat summand of $\mathbb V_{\mathbb C}$ containing the period line
has the form $L\otimes\mathbb W$, where $L$ is fixed and $\mathbb W$
is a rank-two complex local system. In this tensor-product decomposition, one has
\[
  F^2\mathbb V_{\mathbb C}=L^{1,0}\otimes\mathbb W^{1,0},
  \qquad \dim L^{1,0}=1.
\]
In a local trivialization of $\mathbb W$ with a reference fiber $W$, the
period line therefore lies in the fixed projective line
$\mathbb P(L^{1,0}\otimes W)\simeq\mathbb P^1$. Since the period line
determines the point in the period domain $\mathcal D$, the differential of $f$ has rank
at most one everywhere. This contradicts the hypothesis on $\mathrm{rk}(df)$.

Finally, if $f$ is quasi-finite, its fibers have dimension zero, and thus one has $\operatorname{rk}(df_x)
=\dim X\geq2$ at a general point $x\in X$.
\end{proof}

\begin{remark}
  Without the rank assumption Theorem~\ref{h-deform-rigidity_thm_ball_case} no longer holds true. Here is an example which comes from Deraux's construction (\cite{Deraux11}). For the Deligne--Mostow weights $\mu=\frac{1}{8}(3,3,3,3,3,1)$ and $\nu=\frac{1}{8}(3,3,3,7)$ there is a nonconstant holomorphic map of compact ball-quotient orbifolds
\[
f_0:X_\mu=\mathbb B^3/\Gamma_\mu
   \longrightarrow
C_\nu=\mathbb B^1/\Gamma_\nu.
\]
In fact, \(C_\nu\) occurs as a totally geodesic suborbifold of \(X_\mu\), and \(f_0\) is a retraction onto it. For the target weights \(\nu\), the three orbifold orders are
\[
\left(1-\frac{3}{8}-\frac{3}{8}\right)^{-1}=4,
\]
so \(C_\nu\) is the triangle orbifold of signature \((4,4,4)\). The group \(\Delta(4,4,4)\) is one of the compact arithmetic triangle groups in Takeuchi's classification. Note that the trace field of $\Delta(4,4,4)$ is $F=\mathbb{Q} (\sqrt{2})$, and the Shimura curve $C_\nu$ is constructed from a quaternion algebra \(B/F\), which splits at one real place and ramified at other places. Afer taking some finite orbifold cover $C'\to C_\nu$, the construction of Saito--Zucker in \cite{Saito-Zucker} gives us a non-rigid
VHS of K3 type (they call it the basic exmaple in Sect.~6 of the {\it loc. cit.}). In particular, this gives a {\it non-rigid} map from $C_\nu$ to $Y= \Gamma \backslash \mathcal{D}_T$. Therefore, the composite map from $X_\mu=\Gamma_\mu\backslash\mathbb B^3$ to $Y$ is never rigid.
\end{remark}

\begin{thm}[Intermediate Jacobian case]\label{h-deform-rigidity_thm_intJac_case}
The period map
          \[
   p:\, M_{n,d} \to \mathcal{A}_{N(n,d)}
 \]
 induced by taking the intermediate Jacobian of smooth hypersurfaces listed in Sect.~\ref{period:hypersurface} is $h$-deformation rigid.
\end{thm}

\begin{proof}
Suppose that $p$ is homotopic to a distinct holomorphic orbifold map $q$.
Pull this homotopy back to the smooth parameter space $U_{n,d}$
of Sect.~\ref{period:hypersurface}. Choose a neat principal congruence
subgroup $\Gamma\subset\mathrm{Sp}_{2N(n,d)}(\mathbb Z)$ and the
connected finite \'etale cover $X\to U_{n,d}$ induced by the monodromy.
The homotopy lifts to $Y=\mathfrak H_{N(n,d)}/\Gamma$; write $f,g:X\to Y$
for its endpoints. These maps remain distinct.
According to \cite[\S4]{Nog88}, $X$ admits a complete finite-volume
K\"ahler metric for which $f,g$ have finite energy
\cite[Lemma~4.2]{Nog88}. Thus Theorem~\ref{SchYau} and
\cite[Lemma~4.6]{Nog88} give a one-parameter family of holomorphic maps
$\Psi:X\times[0,1]\to Y$.

At a point where $f$ and $g$ differ, the evaluation of $\Psi$
is a nonconstant geodesic. Hence $f,g$ belong to a positive-dimensional
connected component $Z\subset\mathrm{Hol}(X,Y)$. By
Theorem~\ref{Nog-thm}, $Z$ is smooth and quasi-projective.
As before, we obtain a holomorphic deformation of the period map
\[
 \Phi:X\times Z\longrightarrow Y.
\]
Then we receive an extension of the monodromy representation $\rho$:
\[
 \rho_\Phi:\pi_1(X)\times\pi_1(Z)\longrightarrow G(\mathbb Q),
 \qquad G=\mathrm{Sp}_{2N(n,d)}.
\]
The images of the two fundamental-group factors commute. By Beauville's
monodromy theorem \cite[Th\'eor\`eme~4]{Beauville}, $\rho(\pi_1(X))$
has finite index in $G(\mathbb Z)$ and is Zariski dense in $G$.
Consequently $\rho_\Phi(\pi_1(Z))$ lies in the finite center of $G$.
Since $\rho_\Phi$ takes values in the neat group $\Gamma$, this image is
trivial. The evaluation map $e_x:Z\to Y$ therefore lifts holomorphically
to the bounded realization of Siegel space. Its coordinates are bounded
holomorphic functions on the quasi-projective variety $Z$; they extend
across the boundary of a smooth projective compactification and hence
are constant. This contradicts the evaluation immersion in
Theorem~\ref{Nog-thm} and $\dim Z>0$.
\end{proof}

\section{Finiteness in the equal-rank case}\label{sec:equal-rank-finiteness}

We now combine the $h$-deformation rigidity established in the previous section with the boundedness part of the Parshin--Arakelov strategy in their proof of the Geometric Shafarevich Conjecture.  We write
$\mathrm{Hol}_{\mathrm{nc}}(X,Y)$ for the set of nonconstant holomorphic
maps of complex orbifolds from $X$ to $Y$.

\begin{lem}[Boundedness]\label{lem:finite_type}
Let    $C=\overline C\setminus S$ be a smooth quasi-projective curve,
where $\overline C$ is projective and $S$ is finite. Let   $Y$ be a
locally symmetric variety associated with   a   neat arithmetic group   .
Then holomorphic maps from $C$ to $Y$ are parametrized by a finite-type reduced
scheme   $\mathcal H$    with an evaluation map $C\times\mathcal H\to Y$.
\end{lem}

\begin{proof}
We fix an automorphic line bundle $\lambda$ on $Y$ whose
invariant Hermitian metric has positive curvature.  After replacing
$\lambda$ by a positive tensor power, it is the restriction of an ample
line bundle $\lambda^{\mathrm{BB}}$ on the Baily--Borel compactification
$Y^{\mathrm{BB}}$. Note that the line bundle $\lambda$ is constructed from a fixed integral polarized
variation of Hodge structures (VHS) on $Y$: one uses a faithful Hodge
representation, or the adjoint representation if $Y$ is not of Hodge
type.  Thus $\lambda$ is a fixed tensor product of determinant line
bundles of the Hodge filtration.

For every holomorphic map $h:C\to Y$, the pullback of the VHS on $Y$ is a
polarized VHS on $C$.  Its local
monodromies at the points of $S$ are    quasi-unipotent, and neatness makes
them unipotent  . By \cite[Theorem~3.1]{Peters00}, there are upper bounds for
the degrees of Deligne extensions of Hodge bundles in terms only of their
ranks, the fixed Hodge numbers, and $\deg\Omega^1_{\overline C}(\log S)$.
Since $\lambda$ is a fixed tensor product of their determinant bundles,
there exists a constant $B_0=B_0(\overline C,S,Y,\lambda)$ such that
\begin{equation}\label{eq:ample-degree-bound}
  \deg_{\overline C}\overline h^*\lambda^{\mathrm{BB}}\leq B_0   
\end{equation}
for every such $h$.  Here    $\overline h:\overline C\to Y^{\mathrm{BB}}$   is the
unique extension of $h$, whose existence follows from
the Borel extension theorem for maps to arithmetic quotients    \cite{Borel72}  . Since $\overline C$ and    $Y^{\mathrm{BB}}$
  are projective, $\overline h$ is algebraic by GAGA.    

Recall that the Hom scheme    $\operatorname{Hom}(\overline C,Y^{\mathrm{BB}})$   is locally of
finite type.  With respect to the ample line bundle $\lambda^{\mathrm{BB}}$, the Hilbert
polynomial of the graph of $\overline h$ is determined by $\deg_{\overline C}\overline h^*\lambda^{\mathrm{BB}}$. Then the uniform degree upper bound
\eqref{eq:ample-degree-bound} implies that only finitely many Hilbert polynomials can
occur, and the Hom scheme is
therefore of finite type.  The    condition $\overline h(C)\subset Y$ is constructible: its
complement is the projection of the inverse image of
$Y^{\mathrm{BB}}\setminus Y$ under evaluation on $C$. Decompose this
constructible set into finitely many disjoint locally closed reduced
strata and let   $\mathcal H$    be their disjoint union. Its evaluation on
$C$ factors through $Y$  .
\end{proof}

\begin{thm}\label{thm:equal-rank-finiteness}
Let $X$ and $Y$ be locally symmetric varieties  , with irreducible source lattice  . If
$\mathrm{rk}(X)=\mathrm{rk}(Y)\geq 2$, then
\[
  \#\mathrm{Hol}_{\mathrm{nc}}(X,Y)<\infty .
\]
\end{thm}

\begin{proof}
We first prove the statement when $X$ and $Y$ are connected smooth
locally symmetric varieties associated with {\it neat} arithmetic groups.  Choose a smooth
projective compactification $\overline X$ such that
$\Delta_X:=\overline X-X$ is a simple normal crossing divisor.  By the
Lefschetz theorem on $\pi_1$ of smooth quasi-projective
varieties, a sufficiently general complete-intersection curve
$\overline C\subset \overline X$ (intersecting with $\Delta_X$ transversely) has the
property that, for $C:=\overline C\cap X$ and $S:=\overline C-C$, the natural homomorphism $\pi_1(C)\longrightarrow \pi_1(X)$ is surjective.  

Let    $\mathcal H$ be the finite-type parameter space furnished by
Lemma~\ref{lem:finite_type}, using $Y^{\mathrm{BB}}$ as target compactification  . Restriction to $C$ defines a map
\begin{equation}\label{eq:restriction-to-curve}
  \mathrm{Hol}_{\mathrm{nc}}(X,Y)\longrightarrow \mathcal H(\mathbb C),
  \qquad f\longmapsto \overline{f|_C}.
\end{equation}
This map is injective:  If $f|_C=g|_C$, then $f$ and $g$ have the same
value at a chosen base point of $C$, and the surjectivity of $\pi_1(C) \twoheadrightarrow \pi_1(X)$ implies that they induce the same
homomorphism on $\pi_1(X)$.  The uniqueness theorem of
Borel--Narasimhan \cite{BorelNarasimhan66} then gives $f=g$.

Moreover, we can prove that each connected component of $\mathcal H$ contains the
restriction of at most one map in $\mathrm{Hol}_{\mathrm{nc}}(X,Y)$:
Suppose that $\overline{f|_C}$ and $\overline{g|_C}$ lie in the same
connected component. Since this connected component is path connected, the universal evaluation map
produces a free homotopy between $f|_C$ and $g|_C$.  Their induced
homomorphisms $\pi_1(C)\to\pi_1(Y)$ are therefore conjugate.  By
the surjectivity of $\pi_1(C) \twoheadrightarrow \pi_1(X)$, the homomorphisms
$f_*,g_*:\pi_1(X)\to\pi_1(Y)$ are conjugate as well.  Note that the universal
covers of $X$ and $Y$ are bounded symmetric domains,    so $X$ and $Y$   are
$K(\pi,1)$-spaces. Hence, $f$ and $g$ are freely homotopic.  Since
$\rk(X)=\rk(Y) \geq 2$, the {\it $h$-deformation rigidity theorem} (Theorem~\ref{h-deform-rigidity_thmER}) implies $f=g$. 
As a finite-type scheme has only finitely many connected
components, the theorem follows in the neat case.

Finally, we remove the neatness assumption.  Choose finite orbifold
\'etale covers $X_0\to X$ and $Y_0\to Y$ with $X_0$ and $Y_0$ neat.  For
a map $f:X\to Y$, the pullback
\[
  X_0\times_Y Y_0\longrightarrow X_0
\]
is a finite \'etale cover of degree at most $\deg(Y_0/Y)$.  Since
$\pi_1(X_0)$ is finitely generated, there are only finitely many
isomorphism classes of connected finite \'etale covers of $X_0$ of
bounded degree.  On each connected component $Z$ of this pullback, the
map $f$ lifts to a nonconstant holomorphic map $Z\to Y_0$.  The varieties
$Z$ have the same rank as $X$ and irreducible lattices, so the neat case gives only finitely many
such lifts for each of the finitely many possible covers $Z\to X_0$.
Because a map $f$ is determined by its pullback to the surjective cover
$Z\to X$, only finitely many maps $X\to Y$ can occur.
\end{proof}
It seems that the above approach is hard to make the finiteness to be effective or quantitative, since in the proof one has to choose a sufficiently general curve $C$ in $X$. Consequently, there is no control of the degree of the chosen curve $C$ as well as the cardinality of the connected components of the Hom scheme.

As a direct application of Theorem~\ref{thm:equal-rank-finiteness}, we also have
\begin{cor}\label{cor:finiteness}
$\#\mathrm{Hol}_{\mathrm{nc}}(\mathcal A_2,\mathcal F_4)< \infty$.
\end{cor}

\section{Failure of global rigidity and an alternative approach to the finiteness}\label{sec:failure_alternative}

\subsection{Failure of the global rigidity for the Kummer construction}\label{failure_Kummer}
In this subsection, we show the following. 
\begin{thm}\label{counterexam_kummer}
    There is a non-constant holomorphic map from $\mathcal{A}_2$ to $\mathcal{F}_4$ different from the Kummer construction map as given in \S \ref{Kummer}, modulo $\aut_{orb}(\mathcal{A}_2)$ and $\aut_{orb}(\mathcal{F}_4)$.
\end{thm}


The following lemma is given in \cite[Theorem 3.21]{Milne} for the smooth cases. We give the statement in the category of orbifolds.
\begin{lem}\label{aut_orb}
Let    $\cD$   be a bounded symmetric domain   , let $G=\operatorname{Aut}_{\mathrm{hol}}(\cD)$, and let   $\cD/\Gamma$ be a locally symmetric variety with    $\Gamma\leq G$  . Then  $ \aut_{orb}(\cD/\Gamma)\cong N_G(\Gamma)/\Gamma$, i.e., the orbifold automorphism group is isomorphic to the normalizer of $\Gamma$ in $G$.
\end{lem}
\begin{proof}
Each orbifold automorphism of $\cD/\Gamma$ can be lifted to an automorphism $f\in G$ of $\cD$ and $f\gamma f^{-1}$ is a deck transformation for any $\gamma\in \Gamma$ and therefore lies in $\Gamma$, i.e., $f\in N_G(\Gamma)$. Also, a trivial automorphism is lifted to an element in $\Gamma$. Conversely any $f\in N_G(\Gamma)$ defines an orbifold automorphism of $\cD/\Gamma$. Hence $ \aut_{orb}(\cD/\Gamma)\cong N_G(\Gamma)/\Gamma$.
\end{proof}

Note that for any even lattice $T$ with signature $(2,s)$, $N_{SO_0(2,s)}(\widetilde{\SO_0}(T))\supset \SO_0(T)$. Moreover, from \cite[Corollary 2]{KS22} we know
\begin{pro}\label{normalizer_lem}
    Let $T$ be an even lattice with signature $(2,s)$ containing $U^{2}$ as a primitive sublattice. Then $N_{SO_0(2,s)}(\widetilde{\SO_0}(T))=\SO_0(T)$. 
\end{pro}
Let \[T \cong \langle -4\rangle \oplus U^2 \oplus  E^2_8(-1)\] contains $U^2$ as a primitive sublattice.
  We obtain \[\aut_{orb}(\cD_{T}/\Gamma_0)=\SO_0(T)/\widetilde{\SO_0}(T)\cong  \mathbb{Z}/2\mathbb{Z}.\]

\begin{proof}[Proof of Theorem~\ref{counterexam_kummer}]
Choose a root $a$ of square $-2$ in the first $E_8(-1)$ summand of
$T$. There is a primitive embedding
\[
 P:=U^{\oplus2}\oplus\mathbb Za
       \cong U^{\oplus2}\oplus\langle-2\rangle\hookrightarrow T,
\]
with orthogonal complement
\[
 N\cong\langle-4\rangle\oplus E_7(-1)\oplus E_8(-1).
\]
The primitive exterior-square representation of the integral
symplectic group is
\[
 r:\mathrm{Sp}(4,\mathbb Z)\longrightarrow\SO(P),
\]
which is given by the orthogonal complement of the principal
polarization $\omega$ in $\bigwedge^2\mathbb Z^4$. The exceptional isomorphism
$\mathrm{Sp}(4,\mathbb R)/\{\pm1\}\cong\mathrm{SO}_0(2,3)$ identifies
$\mathbb H_2$ holomorphically with $\mathcal D_P$  .

The group $A_P=P^\vee/P$ has order two and $\operatorname{O}(A_P)$ is trivial. Consequently, under the natural homomorphism $O(P\oplus N)\rightarrow O(A_P)\oplus O(A_N)$, $r(\gamma)\oplus\mathrm{id}_N$
preserves the gluing subgroup defining $T$ as an overlattice of
$P\oplus N$ (see Proposition \ref{orthogonal_disc_iso}). It extends to an integral isometry of $T$, acts trivially
on $A_{T}\subset A_P\oplus A_N$, and hence belongs to $\Gamma_0$. This defines a group homomorphism \[ \mathrm{Sp}(4,\mathbb Z) \to \Gamma_0.\] We obtain an equivariant
holomorphic inclusion $\mathbb H_2\cong\mathcal D_P\hookrightarrow
\mathcal D_{T}$, and hence a nonconstant map
\[
 p_{\mathrm{st}}:\mathcal A_2\longrightarrow\mathcal F_4.
\]

For a map with monodromy $\rho$, let $H$ be the identity component of
the Zariski closure of its image and define its moving lattice by
\[
 S_\rho:=T\cap\bigl((T_{\mathbb Q})^H\bigr)^\perp.
\]
For $p_{\mathrm{st}}$ this is $P$, of discriminant $2$. For the Kummer
map, the orthogonal complement of the Kummer lattice in K3 cohomology
is $U^{\oplus3}(2)$ (see \cite[Chapter~14, Proposition~3.14]{Huybrechts2016}), and the primitive polarization gives $\langle 4\rangle$.
Thus its moving lattice is
\[
 S_K\cong U^{\oplus2}(2)\oplus\langle-4\rangle,
 \qquad |\operatorname{disc}(S_K)|=64.
\]
The monodromy on the sixteen exceptional classes is finite, so it
does not change this description using $H$.

   A composition by an orbifold automorphism of $\mathcal A_2$ leaves
the monodromy image unchanged. By Proposition~\ref{normalizer_lem},
a composition by an automorphism of $\mathcal F_4$ changes the moving
lattice by an integral isometry. Its discriminant is therefore
invariant under both operations. Since discriminant of the moving lattices are not the same, the two maps are
inequivalent  .
\end{proof}



\subsection{Alternative approach to the finiteness of the holomorphic maps} \label{alternative_finite}
In this subsection we give an alternative proof of  $\#\mathrm{Hol}_{\mathrm{nc}}(\mathcal A_2,\mathcal F_4)< \infty$  proved in Corollary \ref{cor:finiteness}. The proof consists of two ingredients. First, Margulis superrigidity of higher-rank lattices gives finiteness of the complex
monodromy types. Second, the Borel--Harish-Chandra theorem gives
finiteness of integral orthogonal conjugacy classes within each
complex type. Together with the $h$-deformation rigidity established
in Theorem~3.6, these results imply the desired finiteness.

Recall that \[
T:=U^{2}\oplus E^2_8(-1)\oplus \langle -4\rangle,
\]
and 
$
\widetilde{\SO_0}(T)\subseteq \SO_0(T)
$
is the discriminant kernel. The natural morphism
\[
\mathcal{F}_4=\cD_{T}/\widetilde{\SO_0}(T)
\longrightarrow
\cD_{T}/\SO_0(T)
\]
is a finite orbifold covering. Since the source $\cA_2$ is connected, every map to $\cD_{T}/\SO_0(T)$ has at most finitely many lifts to $\mathcal{F}_4$. It is therefore enough to prove the finiteness of $\mathrm{Hol}_{\mathrm{nc}}(\cA_2,\cD_{T}/\SO_0(T))$.

For $f\in \mathrm{Hol}_{\mathrm{nc}}(\cA_2,\cD_{T}/\SO_0(T))$, let
\[
\rho_f:\Sp(4,\bZ) \longrightarrow \SO_0(T)
\]
be the induced monodromy representation. It is defined up to $\SO_0(T)$-conjugacy. Both the source and the target are orbifold $K(\pi,1)$-spaces, since their universal orbifold covers are bounded symmetric domains. Thus conjugate monodromy representations determine the same free orbifold homotopy class. By the $h$-deformation rigidity proved in \ref{h-deform-rigidity_thmER}, two maps in $\mathrm{Hol}_{\mathrm{nc}}(\cA_2,\cD_{T}/\SO_0(T))$ inducing conjugate monodromy representations are equal. Hence the monodromy assignment gives an injection
\begin{equation}\label{eq:alternative-monodromy-injection}
\mathrm{Hol}_{\mathrm{nc}}(\cA_2,\cD_{T}/\SO_0(T)) \hookrightarrow \mathscr R/\SO_0(T),
\end{equation}
where
\[
\mathscr R:=\{\rho_f\mid f\in \mathrm{Hol}_{\mathrm{nc}}(\cA_2,\cD_{T}/\SO_0(T))\}
\]
and $\SO_0(T)$ acts by conjugation.

The pullback by $f$ of the standard weight-$2$ variation of Hodge structure  on $\cD_{T}$ is a polarizable rational variation of Hodge structures. Consequently, Deligne's semisimplicity theorem for polarizable variations of Hodge structures \cite[Proposition~1.13]{Deligne87} implies that
\[
\rho_{f,\bC}:\Sp(4,\bZ) \longrightarrow \GL(T_{\mathbb{C}})
\]
is semi-simple for every $f\in \mathrm{Hol}_{\mathrm{nc}}(\cA_2,\cD_{T}/\SO_0(T))$.

We now prove finiteness of the complex monodromy types using
Margulis superrigidity rigidity of higher-rank lattices. 

\begin{pro}\label{prop:finite-complex-monodromy}
Only finitely many $\operatorname{SO}(T_{\mathbb C})$-conjugacy
classes occur among the complexifications of the representations
in $\mathscr R$.
\end{pro}

\begin{proof}
The group $\Sp(4,\bZ)$ is an irreducible lattice in
$\operatorname{Sp}(4,\mathbb R)$, whose real rank is two.
The Margulis superrigidity
implies that $\Sp(4,\bZ)$ has only finitely many isomorphism classes
of irreducible complex representations of each fixed dimension.
There are therefore only finitely many isomorphism classes of
semisimple complex representations of $\Sp(4,\bZ)$ of dimension $21$.

It remains to compare linear and orthogonal conjugacy.
Let $\rho_1,\rho_2:\Sp(4,\bZ)\to\operatorname{SO}(T_{\mathbb C})$
be linearly conjugate, and choose
$g\in\operatorname{GL}(T_{\mathbb C})$ such that
\[
g\rho_1(\delta)g^{-1}=\rho_2(\delta)
\qquad\text{for every }\delta\in\Sp(4,\bZ).
\]

Let $A=g^T g$, which commutes
with $\rho_1(\Sp(4,\bZ))$. Choose a square root $s$ of $A$ that is a polynomial in $A$.
Then $s^T=s$, and $s$ commutes with $\rho_1(\Sp(4,\bZ))$.
Hence $h=gs^{-1}$ still intertwines the representations, and
\[
h^T h=s^{-1}As^{-1}=I.
\]
Thus $h\in\operatorname{O}(T_{\mathbb C})$. Since $\dim T_{\mathbb Q}=21$
is odd, replacing $h$ by $-h$ if necessary gives
$h\in\operatorname{SO}(T_{\mathbb C})$ without changing the
conjugacy. Consequently, linear conjugacy implies
$\operatorname{SO}(T_{\mathbb C})$-conjugacy.
\end{proof}



We recall the following finiteness theorem of Borel and Harish-Chandra.
\begin{thm}[Theorem 6.9 in \cite{BHC62}]\label{BHC_finite_thm}
    Let $G$ be a reductive algebraic group defined over $\mathbb{Q}$ and $\pi: G\rightarrow \GL(W)$  a representation defined over $\mathbb{Q}$. Let $G_\bZ$ be the group of integral points. Suppose $L$ is a lattice in
$W_\mathbb{Q}$ invariant under $G_\mathbb{Z}$ and    $X=G_{\mathbb C}\cdot x$ a Zariski-closed orbit through $x\in W_{\mathbb Q}$  . Then $X\cap L$ consists
of a finite number of orbits of $G_\mathbb{Z}$.
\end{thm}






We now use Theorem \ref{BHC_finite_thm} to show the following proposition.

\begin{thm}\label{thm:alternative-finiteness}
The set $\mathrm{Hol}_{\mathrm{nc}}(\cA_2,\cF_4)$ of nonconstant holomorphic maps 
is finite.
\end{thm}

\begin{proof}
Let $\delta_1,\cdots,\delta_m$ be the generators of $\Sp(4,\bZ)$. For each $\rho \in \mathscr R$, we can associate it to the tuple \[x_\rho=(\rho(\delta_1),\cdots,\rho(\delta_m))\in \operatorname{End}_\mathbb{Z}(T)^m\]
We know $\SO(T_{\mathbb Q})$ acts by simultaneous conjugation on $\operatorname{End}_\mathbb{Q}(T_{\mathbb Q})^m$. From \cite[Theorem 3.6]{Ric88}, $\SO(T_{\mathbb C}).x_\rho$ is Zariski-closed. 
Then from Proposition \ref{prop:finite-complex-monodromy} and Theorem \ref{BHC_finite_thm}, we know $\mathscr R/\SO_0(T)$
is finite. The injection \eqref{eq:alternative-monodromy-injection} shows that $\mathrm{Hol}_{\mathrm{nc}}(\cA_2,\cD_{T}/\SO_0(T))$ is finite. Finally, the finite orbifold covering
\[
\cD_{T}/\widetilde{\SO_0}(T)
\longrightarrow
\cD_{T}/\SO_0(T)
\]
has only finitely many lifts of each map from the connected orbifold $\cA_2$, so $\mathrm{Hol}_{\mathrm{nc}}(\cA_2,\cF_4)$ is finite as well.
\end{proof}

\section{An explicit upper bound}\label{sec:explicit-bound}

In this section we make the finiteness of
$\mathrm{Hol}_{\mathrm{nc}}(\mathcal A_2,\mathcal F_4)$ effective and
prove Theorem~\ref{mainthm}(2).

\begin{thm}\label{thm:explicit-kummer-bound}
One has $\#\mathrm{Hol}_{\mathrm{nc}}(\mathcal A_2,\mathcal F_4)<2^{20249}$.
\end{thm}

We retain $T$, $\Gamma_0=\widetilde{\SO_0}(T)$ and the
set $\mathscr R$ of monodromy representations from
\S\ref{alternative_finite}. In particular,
$\mathcal F_4=\cD_{T}/\Gamma_0$, and each $\rho\in\mathscr R$
is a representation $\Sp(4,\bZ)\to\SO_0(T)$.
We write $\SO(T)$ for the group of determinant-one integral
isometries of $T$.

The proof counts integral monodromy representations. By
$h$-deformation rigidity, a map
$f\in\mathrm{Hol}_{\mathrm{nc}}(\mathcal A_2,\mathcal F_4)$ is
determined by its monodromy up to $\Gamma_0$-conjugacy.
We show that each such representation $\rho$ preserves a finite-index
orthogonal sum $S_\rho\oplus N_\rho\subseteq T$. The first lattice
has signature $(2,3)$ and carries a twist of the primitive
exterior-square action of $\Sp(4,\bZ)$; the second is negative definite
and has finite monodromy. We call them the moving and definite parts.
Their equivariant isometry classes, together with the glue group
$H_\rho=T/(S_\rho\oplus N_\rho)$, determine $\rho$ up to
$\SO(T)$-conjugacy (Lemma~\ref{lem:explicit-reconstruction}).
The integer $c$ constructed in Proposition~\ref{pro:explicit-index}
satisfies $c\,T\subseteq S_\rho\oplus N_\rho$ for every $\rho$,
bounding the discriminants and the order of the glue group.
Table~\ref{tab:explicit-count} gives the resulting count,
$2^{20249}$.

\begin{table}[htbp]
\centering
\renewcommand{\arraystretch}{1.25}
\caption{The factors of the bound in Theorem~\ref{thm:explicit-kummer-bound}}
\label{tab:explicit-count}
\begin{tabular}{|l|l|l|}
\hline
Datum & Bound & Reference\\
\hline
moving part $(S_\rho,\sigma_\rho)$, up to isomorphism
 & $<2^{312}$ & Proposition~\ref{pro:explicit-moving-count}\\
definite part $(N_\rho,\kappa_\rho)$, up to isomorphism
 & $<2^{14925}$ & Proposition~\ref{pro:explicit-definite-count}\\
glue group $H_\rho$, for given moving and definite parts
 & $<2^{5010}$ & \S\ref{subsec:explicit-conclusion}\\
$\Gamma_0$-conjugacy classes in an $\SO(T)$-conjugacy class
 & $\leq4$ & Corollary~\ref{cor:explicit-reduction}\\
\hline
\end{tabular}
\end{table}

\subsection{The moving part, the definite part and the glue group}
\label{subsec:explicit-structure}
Let $e_1,e_2,f_1,f_2$ be the standard symplectic basis of $\bZ^4$ and
$\omega=e_1\wedge f_1+e_2\wedge f_2$. The group $\Sp(4,\bZ)$ acts on
$\bigwedge^2\bZ^4$. It fixes $\omega$ and preserves the pairing
$(x,y)\mapsto x\wedge y\in\bigwedge^4\bZ^4\cong\bZ$, where the
isomorphism sends $e_1\wedge e_2\wedge f_1\wedge f_2$ to $1$. Let
$M\subset\bigwedge^2\bZ^4$ be the orthogonal complement of $\omega$.
It has the basis
\begin{equation}\label{eq:explicit-basis-M}
 b_1=e_1\wedge e_2,\quad b_2=f_1\wedge f_2,\quad
 b_3=e_1\wedge f_2,\quad b_4=e_2\wedge f_1,\quad
 b_5=e_1\wedge f_1-e_2\wedge f_2,
\end{equation}
and the action of $\Sp(4,\bZ)$ on $M$ is the primitive exterior-square
representation $r:\Sp(4,\bZ)\to\GL(M)$ used in the proof of
Theorem~\ref{counterexam_kummer}. We equip $M$ with the $r$-invariant
form $(x,y)_M=-2\,x\wedge y$. Then
\[
 (b_1,b_2)_M=(b_3,b_4)_M=-2,\qquad(b_5,b_5)_M=-4,
\]
and all other products of basis vectors vanish. Hence
\begin{equation}\label{eq:explicit-M}
 M\cong U(2)^{\oplus2}\oplus\langle-4\rangle,\qquad
 |\det M|=64,\qquad (M,M)_M=2\bZ,
\end{equation}
where $(M,M)_M$ denotes the ideal generated by the values of the form.
The lattice $M$ is isometric to the moving lattice $S_K$ of the Kummer
map found in \S\ref{failure_Kummer}.

By a \emph{$\Sp(4,\bZ)$-lattice} we mean a lattice $L$ together with a
homomorphism $\Sp(4,\bZ)\to\mathrm O(L)$. An isomorphism of
$\Sp(4,\bZ)$-lattices is an equivariant isometry.

\begin{pro}\label{pro:monodromy-decomposition}
Let $\rho\in\mathscr R$. There are a character
$\epsilon:\Sp(4,\bZ)\to\{\pm1\}$, a positive rational number $a$ and an
injective $\bQ$-linear map $\varphi:M_{\bQ}\to T_{\bQ}$ with the
following properties, where $S_{\bQ}=\varphi(M_{\bQ})$ and
$N_{\bQ}=S_{\bQ}^{\perp}$.
\begin{enumerate}[label=\textup{(\roman*)}]
 \item $\rho(\delta)\,\varphi=\epsilon(\delta)\,\varphi\,r(\delta)$
 for every $\delta\in\Sp(4,\bZ)$.
 \item $(\varphi(x),\varphi(y))=a\,(x,y)_M$ for all $x,y\in M_{\bQ}$.
 \item $N_{\bQ}$ is negative definite of dimension $16$, and
 $\rho(\Sp(4,\bZ))$ acts on it through a finite group.
\end{enumerate}
In particular, $T_{\bQ}=S_{\bQ}\perp N_{\bQ}$ is a decomposition into
$\rho(\Sp(4,\bZ))$-invariant subspaces, and $S_{\bQ}$ has signature $(2,3)$.
\end{pro}

\begin{proof}
Let $f$ be a map with monodromy $\rho$ and let
$\tilde f:\mathbb H_2\to\cD_{T}$ be its lift to the universal covers,
so that $\tilde f\circ\delta=\rho(\delta)\circ\tilde f$ for
$\delta\in\Sp(4,\bZ)$. By Theorem~\ref{thm:MR}, $\tilde f$ is a totally
geodesic embedding. By the classification of holomorphic totally
geodesic embeddings \cite[Theorems~2 and~4, \S4.4]{Ihara67} (see also
\cite{Satake65}), such an embedding of
$D^{\mathrm{III}}_2\cong D^{\mathrm{IV}}_3$ into $D^{\mathrm{IV}}_{19}$
is a standard orthogonal inclusion: there is a subspace
$S_{\bR}\subset T_{\bR}$ of signature $(2,3)$ such that
$\tilde f(\mathbb H_2)=\cD_{S_{\bR}}$. Its orthogonal complement
$N_{\bR}$ is negative definite of dimension $16$, and $\rho(\Sp(4,\bZ))$
preserves $S_{\bR}$ and $N_{\bR}$.

Under the $r$-equivariant identification $\mathbb H_2\cong\cD_{M_{\bR}}$
from the proof of Theorem~\ref{counterexam_kummer}, the biholomorphism
$\tilde f:\cD_{M_{\bR}}\to\cD_{S_{\bR}}$ is induced by a linear
isometry $\psi:M_{\bR}\to S_{\bR}$. For $\delta\in\Sp(4,\bZ)$, the
isometries $\rho(\delta)|_{S_{\bR}}$ and $\psi\,r(\delta)\,\psi^{-1}$
induce the same automorphism of $\cD_{S_{\bR}}$, and only $\pm1$ act
trivially on $\cD_{S_{\bR}}$. Hence
\begin{equation}\label{eq:explicit-real-intertwiner}
 \rho(\delta)\,\psi=\epsilon(\delta)\,\psi\,r(\delta),
 \qquad\epsilon(\delta)\in\{\pm1\},
\end{equation}
and $\epsilon$ is a character of $\Sp(4,\bZ)$.

We now replace $\psi$ by a $\mathbb Q$-linear map. Condition~(i) is a system of
linear equations with rational coefficients for
$\varphi\in\operatorname{Hom}(M_{\bQ},T_{\bQ})$, so its real
solutions are spanned by its rational solutions. Let $\varphi'$ be a
real solution. Its component in $N_{\bR}$ intertwines $\epsilon r$ with
a representation whose image is relatively compact. As $\epsilon r$ is
irreducible and its image is not relatively compact, this component
vanishes. Thus $\varphi'$ maps into $S_{\bR}$, and by
\eqref{eq:explicit-real-intertwiner} the endomorphism
$\psi^{-1}\varphi'$ of $M_{\bR}$ commutes with $r(\Sp(4,\bZ))$. Since $r$
is absolutely irreducible (its image is Zariski dense in
$\SO(M_{\bC})$), we get $\varphi'\in\bR\psi$. The real solutions
therefore form the line $\bR\psi$, which contains a rational solution
$\varphi=\lambda\psi$ with $\lambda\in\bR^{\times}$.

This $\varphi$ satisfies (i). As $\psi$ is an isometry, (ii) holds
with $a=\lambda^2$, which is rational because both forms in (ii) are
rational on $M_{\bQ}$. Finally, $S_{\bQ}\otimes\bR=S_{\bR}$ and
$N_{\bQ}\otimes\bR=N_{\bR}$, so $N_{\bQ}$ is negative definite of
dimension $16$, and $\rho(\Sp(4,\bZ))$ preserves the definite lattice
$T\cap N_{\bQ}$, whose isometry group is finite.
\end{proof}

For each $\rho\in\mathscr R$ we fix $\epsilon_\rho$, $a_\rho$ and
$\varphi_\rho$ as in Proposition~\ref{pro:monodromy-decomposition},
with the resulting subspaces $S_{\rho,\bQ}$ and $N_{\rho,\bQ}$, and let
\[
 S_\rho=T\cap S_{\rho,\bQ},\qquad N_\rho=T\cap N_{\rho,\bQ},\qquad
 \sigma_\rho=\rho|_{S_\rho},\qquad\kappa_\rho=\rho|_{N_\rho}.
\]
The lattices $S_\rho$ and $N_\rho$ are primitive in $T$, and each is
the orthogonal complement of the other. We call the $\Sp(4,\bZ)$-lattices
$(S_\rho,\sigma_\rho)$ and $(N_\rho,\kappa_\rho)$ the \emph{moving part}
and the \emph{definite part} of $\rho$.

Orthogonal projection embeds $T$ into $S_\rho^\vee\oplus N_\rho^\vee$.
The \emph{glue group}
\[
 H_\rho=T/(S_\rho\oplus N_\rho)\subseteq A_{S_\rho}\oplus A_{N_\rho}
\]
therefore determines $T$, as its preimage in
$S_\rho^\vee\oplus N_\rho^\vee$, and the projection
$H_\rho\to A_{S_\rho}$ is injective
(Proposition~\ref{orthogonal_disc_iso}). We call
\[
 \mathcal G(\rho)=(S_\rho,\sigma_\rho;\,N_\rho,\kappa_\rho;\,H_\rho)
\]
the \emph{gluing datum} of $\rho$. Two gluing data
$(S,\sigma;N,\kappa;H)$ and $(S',\sigma';N',\kappa';H')$ are
\emph{isomorphic} if there are isomorphisms of $\Sp(4,\bZ)$-lattices
$\mu:(S,\sigma)\to(S',\sigma')$ and $\eta:(N,\kappa)\to(N',\kappa')$ such
that the induced isomorphism $A_S\oplus A_N\to A_{S'}\oplus A_{N'}$
maps $H$ onto $H'$.

\begin{lem}\label{lem:explicit-reconstruction}
Let $\rho,\rho'\in\mathscr R$. If the gluing data $\mathcal G(\rho)$
and $\mathcal G(\rho')$ are isomorphic, then $\rho'=g\rho g^{-1}$ for
some $g\in\SO(T)$.
\end{lem}

\begin{proof}
Let $\mu$ and $\eta$ be as in the definition, and let $g=\mu\oplus\eta$,
extended $\bQ$-linearly to an isometry of
$T_{\bQ}=S_{\rho,\bQ}\perp N_{\rho,\bQ}
=S_{\rho',\bQ}\perp N_{\rho',\bQ}$.
It maps $S_\rho^\vee\oplus N_\rho^\vee$ onto
$S_{\rho'}^\vee\oplus N_{\rho'}^\vee$, and the induced map on
discriminant groups carries $H_\rho$ onto $H_{\rho'}$. Since $T$ is
the preimage of $H_\rho$ in the first group and of $H_{\rho'}$ in the
second, $g(T)=T$. The equivariance of $\mu$ and
$\eta$ gives $g\rho(\delta)=\rho'(\delta)g$ for every $\delta\in\Sp(4,\bZ)$.
As $T$ has odd rank, one of $\pm g$ has determinant one, and both
define the same conjugation.
\end{proof}

\begin{cor}\label{cor:explicit-reduction}
The number of nonconstant holomorphic maps $\mathcal A_2\to\mathcal F_4$
is at most four times the number of isomorphism classes of gluing data
$\mathcal G(\rho)$, $\rho\in\mathscr R$.
\end{cor}

\begin{proof}
Let $f\in\mathrm{Hol}_{\mathrm{nc}}(\mathcal A_2,\mathcal F_4)$, with
monodromy representation $\rho_f:\Sp(4,\bZ)\to\Gamma_0$, defined up to
conjugation by $\Gamma_0$. Composing $f$ with the covering
$\mathcal F_4\to\cD_{T}/\SO_0(T)$ shows that $\rho_f\in\mathscr R$.
As for \eqref{eq:alternative-monodromy-injection},
Theorem~\ref{h-deform-rigidity_thmER} shows that $f$ is determined by
the $\Gamma_0$-conjugacy class of $\rho_f$. By
Lemma~\ref{lem:explicit-reconstruction}, the representations $\rho_f$
with a given isomorphism class of gluing data lie in one
$\SO(T)$-conjugacy class, which is the union of at most
$[\SO(T):\Gamma_0]$ conjugacy classes under $\Gamma_0$. This index is
at most four: $\SO_0(T)$ has index at most two in $\SO(T)$, and
$\Gamma_0$ is the kernel of the action of $\SO_0(T)$ on
$A_{T}\cong\bZ/4\bZ$, whose only automorphisms are $\pm1$.
\end{proof}

\subsection{The index of $S_\rho\oplus N_\rho$ in $T$}
\label{group_ring_compu}
For $\rho\in\mathscr R$ let $e_\rho\in\operatorname{End}(T_{\bQ})$ be
the orthogonal projection onto $S_{\rho,\bQ}$. We construct an integer
$c$, independent of $\rho$, such that $c\,e_\rho$ preserves $T$. The
idea is to find an element $z$ of the group ring $\bZ[\Sp(4,\bZ)]$ which
acts on $M$ as multiplication by $c$, whose coefficients sum to zero,
and which is a combination of elements of $\Sp(4,\bZ)$ on which every
$\epsilon_\rho$ and every $\kappa_\rho$ is trivial. Then $\rho(z)$ is
$c$ on $S_{\rho,\bQ}$ and $0$ on $N_{\rho,\bQ}$, that is,
$\rho(z)=c\,e_\rho$. The next lemma
provides the required elements of $\Sp(4,\bZ)$: the $Q$th powers.

\begin{lem}\label{lem:explicit-finite-order}
Let
\[
 Q=2^5\cdot3^2\cdot5\cdot7\cdot11\cdot13\cdot17=24504480<2^{25}.
\]
Every element of finite order of $\GL_{16}(\bZ)$ has order dividing $Q$.
\end{lem}

\begin{proof}
Let $A\in\GL_{16}(\bZ)$ have finite order. Then $A$ is diagonalizable
over $\bC$, and its order is the least common multiple of the orders
of its eigenvalues. If a primitive $n$th root of unity is an
eigenvalue, the cyclotomic polynomial $\Phi_n$ divides the
characteristic polynomial of $A$, so $\deg\Phi_n\leq16$. For a prime
power $p^k$ dividing $n$, the degree
$\deg\Phi_{p^k}=p^{k-1}(p-1)$ divides $\deg\Phi_n$, so
$p^{k-1}(p-1)\leq16$. This holds exactly for the divisors $p^k$ of
$2^5$, $3^2$, $5$, $7$, $11$, $13$ and $17$. Hence $n$ divides $Q$.
\end{proof}

For a symmetric integral $2\times2$ matrix $B$ set
\begin{equation}\label{equ:def_u_l}
 u(B)=\begin{pmatrix}I_2&B\\0&I_2\end{pmatrix},\qquad
 \ell(B)=\begin{pmatrix}I_2&0\\B&I_2\end{pmatrix},
\end{equation}
in the basis $e_1,e_2,f_1,f_2$, and let $E_{ij}$ denote the
$2\times2$ matrix units. For an integer $q\geq1$ consider the six
elements of $\Sp(4,\bZ)$
\[
 \begin{aligned}
 g_1(q)&=u(qE_{11}),& g_2(q)&=u(qE_{22}),& g_3(q)&=u(q(E_{12}+E_{21})),\\
 g_4(q)&=\ell(qE_{11}),& g_5(q)&=\ell(qE_{22}),& g_6(q)&=\ell(q(E_{12}+E_{21})).
 \end{aligned}
\]
They satisfy
\begin{equation}\label{eq:explicit-qth-power}
 g_i(q)=g_i(1)^q\qquad(1\leq i\leq6).
\end{equation}
In the basis \eqref{eq:explicit-basis-M} of $M$, the matrices
$R_i(q)=r(g_i(q))$ are
{\small
\[
 \begin{aligned}
 R_1(q)&=\begin{pmatrix*}[r]
 1&0&0&-q&0\\0&1&0&0&0\\0&q&1&0&0\\0&0&0&1&0\\0&0&0&0&1
 \end{pmatrix*},&
 R_2(q)&=\begin{pmatrix*}[r]
 1&0&q&0&0\\0&1&0&0&0\\0&0&1&0&0\\0&-q&0&1&0\\0&0&0&0&1
 \end{pmatrix*},\\[4pt]
 R_3(q)&=\begin{pmatrix*}[r]
 1&-q^2&0&0&2q\\0&1&0&0&0\\0&0&1&0&0\\0&0&0&1&0\\0&-q&0&0&1
 \end{pmatrix*},&
 R_4(q)&=\begin{pmatrix*}[r]
 1&0&0&0&0\\0&1&q&0&0\\0&0&1&0&0\\-q&0&0&1&0\\0&0&0&0&1
 \end{pmatrix*},\\[4pt]
 R_5(q)&=\begin{pmatrix*}[r]
 1&0&0&0&0\\0&1&0&-q&0\\q&0&1&0&0\\0&0&0&1&0\\0&0&0&0&1
 \end{pmatrix*},&
 R_6(q)&=\begin{pmatrix*}[r]
 1&0&0&0&0\\-q^2&1&0&0&-2q\\0&0&1&0&0\\0&0&0&1&0\\q&0&0&0&1
 \end{pmatrix*}.
 \end{aligned}
\]
}%
Let $\mathcal P(q)\subset\operatorname{End}_{\bZ}(M)$ be the set of
matrices $r(w)$, where $w$ runs through the products of at most three
of the twelve elements $g_i(q)^{\pm1}$. We identify
$\operatorname{End}_{\bZ}(M)$ with $\bZ^{25}$ by means of the basis
\eqref{eq:explicit-basis-M}.

We shall use Smith normal forms in the following way. An integral
$m\times n$ matrix $A$ of rank $m$ has invariant factors
$d_1\mid d_2\mid\cdots\mid d_m$: there are $U\in\GL_m(\bZ)$ and
$V\in\GL_n(\bZ)$ with $UAV=(\diag(d_1,\ldots,d_m)\mid0)$
\cite[Chapter~12]{DF04}. The lattice spanned by the columns of $A$ is
then $U^{-1}(d_1\bZ\oplus\cdots\oplus d_m\bZ)$, so that
\begin{equation}\label{eq:explicit-snf}
 d_m\bZ^m\subseteq A\bZ^n.
\end{equation}
The two Smith normal forms in Lemmas~\ref{lem:explicit-integral-order}
and~\ref{lem:explicit-level-two} were computed in exact integer
arithmetic by the accompanying script.

\begin{lem}\label{lem:explicit-integral-order}
Let $P_1,\ldots,P_n$ be the elements of $\mathcal P(Q)$. There are
integers $z_1,\ldots,z_n$ such that
\[
 \sum_{j=1}^n z_jP_j=Q^4\,\id_M,\qquad\sum_{j=1}^n z_j=0.
\]
\end{lem}

\begin{proof}
Let $A$ be the integral matrix with $26$ rows whose $j$th column
consists of the $25$ entries of $P_j$ followed by the entry $1$. The
computation gives $n=877$, and $A$ has rank $26$ with invariant factors
\[
 1,\quad Q^{[6]},\quad(Q^2)^{[9]},\quad(2Q^2)^{[3]},\quad
 (Q^3)^{[4]},\quad(2Q^3)^{[2]},\quad Q^4,
\]
where $x^{[j]}$ denotes $j$ repetitions of $x$. By
\eqref{eq:explicit-snf}, the vector $(Q^4\,\id_M,0)\in Q^4\bZ^{26}$ is
an integral combination of the columns of $A$.
\end{proof}

\begin{pro}\label{pro:explicit-index}
Let $c=Q^4<2^{100}$. For every $\rho\in\mathscr R$ the endomorphism
$c\,e_\rho$ preserves $T$. Consequently
$c\,T\subseteq S_\rho\oplus N_\rho$, the glue group satisfies
$|H_\rho|\leq c^5$, and
\begin{equation}\label{eq:explicit-determinant-product}
 |A_{S_\rho}|\,|A_{N_\rho}|=4\,|H_\rho|^2\leq4c^{10}<2^{1002}.
\end{equation}
\end{pro}

\begin{proof}
Let $z_j$ and $P_j=r(w_j)$ be as in
Lemma~\ref{lem:explicit-integral-order}, where each $w_j$ is a product
of elements $g_i(Q)^{\pm1}$, and put
$z=\sum_jz_jw_j\in\bZ[\Sp(4,\bZ)]$. By \eqref{eq:explicit-qth-power}, each
$g_i(Q)$ is a $Q$th power in $\Sp(4,\bZ)$. Hence $\epsilon_\rho(g_i(Q))=1$,
because $Q$ is even, and $\kappa_\rho(g_i(Q))=1$, by
Lemma~\ref{lem:explicit-finite-order} applied to the element
$\kappa_\rho(g_i(1))$ of finite order in
$\GL(N_\rho)\cong\GL_{16}(\bZ)$. Consequently $\rho(w_j)$ is the
identity on $N_{\rho,\bQ}$, and
$\rho(w_j)\,\varphi_\rho=\varphi_\rho\,P_j$ by
Proposition~\ref{pro:monodromy-decomposition}(i). Therefore
$\rho(z)=\sum_jz_j\rho(w_j)$ acts on $N_{\rho,\bQ}$ as multiplication
by $\sum_jz_j=0$ and on $S_{\rho,\bQ}$ as
$\varphi_\rho\bigl(\sum_jz_jP_j\bigr)\varphi_\rho^{-1}=c$. Thus
$c\,e_\rho=\rho(z)$, which lies in
$\bZ[\rho(\Sp(4,\bZ))]\subseteq\operatorname{End}_{\bZ}(T)$.

For $t\in T$ we obtain $c\,e_\rho(t)\in T\cap S_{\rho,\bQ}=S_\rho$
and $c\,(1-e_\rho)(t)\in T\cap N_{\rho,\bQ}=N_\rho$, hence
$ct\in S_\rho\oplus N_\rho$. So $c$ annihilates $H_\rho$. Since
$H_\rho$ embeds into $A_{S_\rho}$, a group with five generators, it is
generated by at most five elements of order dividing $c$, and
$|H_\rho|\leq c^5$. Finally, Proposition~\ref{finite_index_disc} for
$S_\rho\oplus N_\rho\subseteq T$ gives
$|A_{S_\rho}|\,|A_{N_\rho}|=|A_{T}|\,|H_\rho|^2=4\,|H_\rho|^2$.
\end{proof}

\subsection{Counting the moving parts}\label{moving_lattice_counting}
Multiplying $\varphi_\rho$ by a nonzero rational number $\lambda$
preserves the properties in
Proposition~\ref{pro:monodromy-decomposition}, with $a_\rho$ replaced
by $\lambda^2a_\rho$, and replaces the lattice
$L_\rho=\varphi_\rho^{-1}(S_\rho)\subset M_{\bQ}$ by
$\lambda^{-1}L_\rho$. We therefore normalize $\varphi_\rho$ so that
\begin{equation}\label{eq:explicit-normalization}
 L_\rho\subseteq M\qquad\text{and}\qquad
 L_\rho\not\subseteq pM\ \text{ for every prime }p.
\end{equation}
By Proposition~\ref{pro:monodromy-decomposition}(i) and~(ii),
$\varphi_\rho$ restricts to an isomorphism of $\Sp(4,\bZ)$-lattices
\begin{equation}\label{eq:explicit-normal-form}
 \bigl(L_\rho,\ \epsilon_\rho\,r|_{L_\rho},\ a_\rho(\cdot,\cdot)_M\bigr)
 \xrightarrow{\ \sim\ }(S_\rho,\sigma_\rho).
\end{equation}
The moving part of $\rho$ is thus determined, up to isomorphism, by the
triple $(\epsilon_\rho,L_\rho,a_\rho)$, and we count these triples.
The lattice $L_\rho$ is controlled by the following computation.

\begin{lem}\label{lem:explicit-level-two}
The lattice spanned by $\mathcal P(2)$ in $\operatorname{End}_{\bZ}(M)$
contains $16\operatorname{End}_{\bZ}(M)$.
\end{lem}

\begin{proof}
The set $\mathcal P(2)$ has $877$ elements. The integral
$25\times877$ matrix whose columns are these elements has rank $25$ and
invariant factors
\[
 1,\quad2^{[6]},\quad4^{[9]},\quad8^{[7]},\quad16^{[2]}.
\]
The claim follows from \eqref{eq:explicit-snf}.
\end{proof}

\begin{pro}\label{pro:explicit-moving-count}
The moving parts $(S_\rho,\sigma_\rho)$, $\rho\in\mathscr R$, belong to
fewer than $2^{312}$ isomorphism classes of $\Sp(4,\bZ)$-lattices.
\end{pro}

\begin{proof}
By \eqref{eq:explicit-normal-form} it suffices to bound the number of
triples $(\epsilon_\rho,L_\rho,a_\rho)$.

\emph{The lattice $L_\rho$.} Since $S_\rho$ is $\rho(\Sp(4,\bZ))$-stable,
$L_\rho$ is stable under $\epsilon_\rho(\delta)r(\delta)$, hence under
$r(\delta)$, for every $\delta\in\Sp(4,\bZ)$. By
Lemma~\ref{lem:explicit-level-two}, it is therefore stable under
$16\operatorname{End}_{\bZ}(M)$. By the elementary divisor theorem
there are a basis $m_1,\ldots,m_5$ of $M$ and positive integers
$d_1\mid\cdots\mid d_5$ with $L_\rho=\bigoplus_i\bZ d_im_i$, and
\eqref{eq:explicit-normalization} forces $d_1=1$, so $m_1\in L_\rho$.
Given $m\in M$, let $J\in\operatorname{End}_{\bZ}(M)$ be defined by
$J(m_1)=m$ and $J(m_i)=0$ for $i\geq2$. Then
$16m=(16J)(m_1)\in L_\rho$. Therefore
\begin{equation}\label{eq:explicit-moving-sandwich}
 16M\subseteq L_\rho\subseteq M.
\end{equation}
Such a lattice corresponds to a subgroup of
$M/16M\cong(\bZ/16\bZ)^5$, and every such subgroup is generated by at
most five elements. There are $16^{25}=2^{100}$ quintuples of elements
of $(\bZ/16\bZ)^5$, hence at most $2^{100}$ lattices $L_\rho$.

\emph{The scalar $a_\rho$.} Write $a_\rho=\alpha/\beta$ in lowest
terms. The form $a_\rho(\cdot,\cdot)_M$ is integral on
$L_\rho\supseteq16M$, and by \eqref{eq:explicit-M} its values on $16M$
generate $512\,a_\rho\bZ$. Hence $512\,a_\rho\in\bZ$, that is,
$\beta$ divides $512$. Next,
\eqref{eq:explicit-normal-form}, Proposition~\ref{finite_index_disc}
and \eqref{eq:explicit-M} give
\[
 |A_{S_\rho}|=a_\rho^5\,[M:L_\rho]^2\,|\det M|
 =64\,a_\rho^5\,[M:L_\rho]^2.
\]
Since $|A_{S_\rho}|<2^{1002}$ by
\eqref{eq:explicit-determinant-product}, we get $a_\rho^5<2^{996}$ and
$a_\rho<2^{200}$. For each of the ten divisors $\beta$ of $512$ there
are fewer than $2^{200}\beta$ numerators $\alpha$, so the number of
possible values of $a_\rho$ is less than
\[
 2^{200}(1+2+\cdots+512)<2^{210}.
\]

\emph{The character $\epsilon_\rho$.} The group $\Sp(4,\bZ)$ is generated
by two elements \cite{Bender80}, so it has at most four characters
with values in $\{\pm1\}$.

Altogether there are fewer than $4\cdot2^{100}\cdot2^{210}=2^{312}$
triples.
\end{proof}

\subsection{Counting the definite parts}\label{fix_lattice_count}
We first count the lattices $N_\rho$, using Minkowski's second
theorem, and then the actions $\kappa_\rho$ on each of them.

\begin{thm}[Minkowski's second theorem
{\cite[Theorem~23.1]{Gru07}}]\label{minko_second}
Let $\Lambda\subset\bR^n$ be a lattice of rank $n$, and let
$K\subset\bR^n$ be a compact convex set with nonempty interior,
symmetric about the origin. Its successive minima are
\[
 \lambda_i=\inf\bigl\{t>0\mid tK\text{ contains }i
 \text{ linearly independent points of }\Lambda\bigr\}
 \qquad(1\leq i\leq n).
\]
Then
\[
 \lambda_1\cdots\lambda_n\vol(K)\leq2^n\vol(\bR^n/\Lambda),\qquad
 \vol(\bR^n/\Lambda)=\sqrt{|\operatorname{disc}(\Lambda)|}.
\]
\end{thm}

\begin{lem}\label{lem:finite-orthogonal-group}
Let $N$ be a definite lattice of rank $16$. Then
$|\mathrm O(N)|<3^{256}<2^{406}$.
\end{lem}

\begin{proof}
The group $\mathrm O(N)$ is a finite subgroup of
$\GL(N)\cong\GL_{16}(\bZ)$. We show that reduction modulo $3$ is
injective on finite subgroups of $\GL_{16}(\bZ)$; then
$|\mathrm O(N)|\leq|\GL_{16}(\bF_3)|<3^{256}$, and $3^{256}<2^{406}$
because $256\log_23<405.8$.

Suppose that $A\neq I$ has finite order and $A\equiv I\pmod3$.
Replacing $A$ by a power, we may assume that its order is a prime $p$.
Write $A=I+3^kB$ with $k\geq1$ and $B\not\equiv0\pmod3$. Expanding
$A^p=I$ gives
\[
 p\,3^kB=-\sum_{i=2}^{p}\binom pi3^{ik}B^i.
\]
If $p\neq3$, the right-hand side is divisible by $3^{k+1}$ and the
left-hand side is not. If $p=3$, the right-hand side
$-3^{2k+1}B^2-3^{3k}B^3$ is divisible by $3^{k+2}$ and the left-hand
side $3^{k+1}B$ is not.
\end{proof}

\begin{pro}\label{pro:explicit-definite-count}
The definite parts $(N_\rho,\kappa_\rho)$, $\rho\in\mathscr R$, belong
to fewer than $2^{14925}$ isomorphism classes of $\Sp(4,\bZ)$-lattices.
\end{pro}

\begin{proof}
Let $\rho\in\mathscr R$. The lattice $\Lambda=N_\rho(-1)$ is positive
definite of rank $16$, and
$D:=\det\Lambda=|A_{N_\rho}|<2^{1002}$ by
\eqref{eq:explicit-determinant-product}.

\emph{Step 1: a sublattice with short basis.} We apply
Theorem~\ref{minko_second} to $\Lambda\subset\Lambda_{\bR}\cong\bR^{16}$
and the unit ball $K$, whose volume is $\pi^8/8!>2^{-3}$. Choose
linearly independent vectors $v_1,\ldots,v_{16}\in\Lambda$ with
$\|v_i\|=\lambda_i$, and put $\nu_i=(v_i,v_i)$. Then
\begin{equation}\label{eq:explicit-minkowski}
 \prod_i\|v_i\|\leq\frac{2^{16}\sqrt D}{\vol(K)}<2^{19}\sqrt D,
 \qquad
 \prod_i\nu_i<2^{38}D<2^{1040}=:X.
\end{equation}
Let $\Lambda'=\bigoplus_i\bZ v_i\subseteq\Lambda$. The lattice
$\Lambda$ is determined, up to isometry, by the Gram matrix
$\bigl((v_i,v_j)\bigr)$ of $\Lambda'$ and by the subgroup
$\Lambda/\Lambda'$ of $A_{\Lambda'}$. We count these two objects in
turn.

\emph{Step 2: the Gram matrix of $\Lambda'$.} The diagonal entries
$\nu_1,\ldots,\nu_{16}$ are positive integers with product less than
$X$. Given $\nu_1,\ldots,\nu_{15}$, there are at most
$X/(\nu_1\cdots\nu_{15})$ values for $\nu_{16}$, so the number of possible diagonals is at most
\[
 \sum_{\nu_1=1}^{X}\cdots\sum_{\nu_{15}=1}^{X}
 \frac{X}{\nu_1\cdots\nu_{15}}
 =X\Bigl(\sum_{j=1}^{X}\frac1j\Bigr)^{15}
 \leq X(1+\log X)^{15}<2^{1040}\cdot(2^{10})^{15}=2^{1190},
\]
where we used $1+\log X<722<2^{10}$. For a given diagonal, the entry
$(v_i,v_j)$ with $i<j$ is an integer of absolute value at most
$\sqrt{\nu_i\nu_j}$ by the Cauchy--Schwarz inequality, which leaves at
most $2\sqrt{\nu_i\nu_j}+1\leq3\sqrt{\nu_i\nu_j}$ values. The number of
choices for the $120$ entries above the diagonal is therefore at most
\[
 \prod_{i<j}3\sqrt{\nu_i\nu_j}=3^{120}\Bigl(\prod_i\nu_i\Bigr)^{15/2}
 <2^{240}\cdot2^{7800}=2^{8040}.
\]
Hence fewer than $2^{1190+8040}=2^{9230}$ Gram matrices occur.

\emph{Step 3: from $\Lambda'$ to $\Lambda$.} By
Proposition~\ref{finite_index_disc}, Hadamard's inequality
$\det\Lambda'\leq\prod_i\nu_i$ and \eqref{eq:explicit-minkowski},
\[
 [\Lambda:\Lambda']=\sqrt{\frac{\det\Lambda'}{D}}
 \leq\frac{\prod_i\|v_i\|}{\sqrt D}<2^{19}.
\]
The group $A_{\Lambda'}$ is generated by sixteen elements. Hence its
$k$-torsion subgroup has at most $k^{16}$ elements, and each of its
subgroups is generated by at most sixteen elements. A subgroup of
order $k$ lies in the $k$-torsion subgroup, so there are at most
$(k^{16})^{16}=k^{256}$ of them. As $\Lambda/\Lambda'$ has order
$[\Lambda:\Lambda']<2^{19}$, the number of lattices $\Lambda$ for a
given Gram matrix is less than
\[
 \sum_{1\leq k<2^{19}}k^{256}<2^{19}\cdot2^{19\cdot256}=2^{4883}.
\]
Thus the lattices $N_\rho$ belong to fewer than
$2^{9230+4883}=2^{14113}$ isometry classes.

\emph{Step 4: the action.} On a given lattice $N$, a homomorphism
$\kappa:\Sp(4,\bZ)\to\mathrm O(N)$ is determined by the images of two
generators of $\Sp(4,\bZ)$. By Lemma~\ref{lem:finite-orthogonal-group}
there are fewer than $(2^{406})^2=2^{812}$ such homomorphisms. The
number of isomorphism classes of definite parts is therefore less than
$2^{14113+812}=2^{14925}$.
\end{proof}

\subsection{The glue groups and the proof of
Theorem~\ref{thm:explicit-kummer-bound}}
\label{subsec:explicit-conclusion}
\begin{proof}[Proof of Theorem~\ref{thm:explicit-kummer-bound}]
Fix a $\Sp(4,\bZ)$-lattice $(S,\sigma)$ in one of the isomorphism classes
of Proposition~\ref{pro:explicit-moving-count} and a $\Sp(4,\bZ)$-lattice
$(N,\kappa)$ in one of the classes of
Proposition~\ref{pro:explicit-definite-count}. Every gluing datum
$\mathcal G(\rho)$ whose moving and definite parts are isomorphic to
$(S,\sigma)$ and $(N,\kappa)$ is isomorphic to
$(S,\sigma;N,\kappa;H)$ for a subgroup $H\subseteq A_S\oplus A_N$
which projects injectively to $A_S$, and then
$|A_S|\,|A_N|<2^{1002}$ by \eqref{eq:explicit-determinant-product}.
Such a subgroup is generated by at most five elements, so the number
of possibilities for $H$ is at most
\[
 \bigl(|A_S|\,|A_N|\bigr)^5<2^{5010}.
\]
Together with Propositions~\ref{pro:explicit-moving-count}
and~\ref{pro:explicit-definite-count}, this shows that the gluing data
$\mathcal G(\rho)$, $\rho\in\mathscr R$, belong to fewer than
\[
 2^{312}\cdot2^{14925}\cdot2^{5010}=2^{20247}
\]
isomorphism classes. By Corollary~\ref{cor:explicit-reduction},
\[
 \#\mathrm{Hol}_{\mathrm{nc}}(\mathcal A_2,\mathcal F_4)
 <4\cdot2^{20247}=2^{20249}.\qedhere
\]
\end{proof}

\section{Discussions for more examples}\label{sec4}
We discuss some rigidity problems for more example of period maps in this section, which serve as the future research. 

\subsection{The case where global rigidity may be far from being true }
It is well-known that $\pi_1^{orb}(\mathrm{UC}_n)=B_n$ is the classical Braid group. More precisely, $B_n$ is generated by $\sigma_1,\cdots,\sigma_{n-1}$ subjective to the relations 
\begin{equation}
   \sigma_i\sigma_j=\sigma_j\sigma_i \ \hbox{for}\ |i-j|\ge 2,\ \   \sigma_i\sigma_{i+1}\sigma_i=\sigma_{i+1}\sigma_i \sigma_{i+1},\ \hbox{for}\ 1\le i\le n-2.
\end{equation}
Note that the abelianization of $B_n$ is $ B_n^{ab}\cong \bZ$, now we fix an isomorphism.  
Denote $T=U^2\oplus E_8^{2}(-1) \oplus \langle -2l \rangle$ the polarized K3 lattice. Let $g\in \SO_0(T)$ and $m\in \bZ_{>0}$,  given any group nontrivial homomorphism $\epsilon:\ B_n^{ab} \rightarrow \bZ$, define the map
\begin{equation}\label{BGH}
    \rho_{g,m,\epsilon}: \ B_n \rightarrow \SO_0(T),\ x\mapsto g^{m\cdot \epsilon(x)}.
\end{equation}
   Writing $k=\epsilon(\sigma_i)\ne0$, we have $\rho_{g,m,\epsilon}(\sigma_i)=g^{mk}$  . It is easy to check that (\ref{BGH}) gives a group homomorphism. 
    Identify $U^2$ with $\mathbb Z^2\oplus(\mathbb Z^2)^\vee$ with its hyperbolic pairing. Now take $g=A\oplus(A^{-1})^t\oplus\mathrm{Id}$   where 
 \[ A= \begin{pmatrix}
     2 & 1\\
     1 &  1
 \end{pmatrix}\] and then we have the following.
\begin{pro}
   If $m \neq l$, then  $\rho_{g,m,\epsilon}$ and $\rho_{g,l,\epsilon}$ are  not conjugate to each other under the conjugate action of $\SO_0(T)$. In particular, there are infinitely many elements in $\homo(B_n,\SO_0(T))/\SO_0(T)$.
\end{pro}
\begin{proof}
Suppose that    $\rho_{g,m,\epsilon}$ and $\rho_{g,l,\epsilon}$ are
conjugate. Their values at $\sigma_1$ have equal traces. If $\lambda>1$
is the larger eigenvalue of $A$, then
\[
 \operatorname{tr}(g^{mk})=2(\lambda^{m|k|}+\lambda^{-m|k|})
                         +\operatorname{rank}(T)-4,
\]  
   which is strictly increasing with $m>0$. Hence $m=l$  .
\end{proof}

   A standard braid generator is a meridian of the collision divisor
in the configuration space. If   a holomorphic map    to the arithmetic
period quotient induced   $\rho_{g,m,\epsilon}$   , the pullback of its
polarized integral VHS would have quasi-unipotent monodromy around
this divisor. The eigenvalue $\lambda^{mk}$ is not a root of unity.
Thus these homomorphisms cannot   be induced by    holomorphic maps  .

Nevertheless, one may try to construct some other group homomorphisms and show that they are induced by holomorphic maps. Then the  failure of group theoretic rigidity  suggests that there may be infinitely many  non-constant holomorphic maps from $\mathrm{UC}_n$ to the locally symmetric $\cD/\Gamma$ different from the period maps constructed in \S \ref{Periodmap:points}. 

\subsection{Rigidity problems for period maps of  moduli of low genus curves }

Besides constructing the period map $M_g \rightarrow \cA_g$ via Jacobian variety, Kondo and his coauthors have found other interesting  period maps for the moduli space of low genus curves ($3\le g\le 6$).  Let's recall their constructions:

\subsubsection*{$g=3$}
Let $V:=|\cO_{\bP^2}(4)|$ be the space of plane quartic curves and    let \[\widetilde{V}:=\mathrm{Bl}_{\PGL_3\cdot 2Q}|\cO_{\bP^2}(4)|\]   be the Kirwan type blowup along the locus of double smooth quadrics $2Q\in V$. Under suitable polarization, there is an open subset $M_3\subset \widetilde{V} \q \PGL_3$ parametrizing all smooth curves of genus $3$ 
, see  \cite{Kondo00} and \cite[Theorem 3.5]{Michela09}. And by associating a K3 surface with polarization of degree $4$ for curve $C\in M_3$,  there is  a period map 
\begin{equation}\label{genus3}
    p: \  M_3 \rightarrow \cD_4/\Gamma
\end{equation}
and $p$ factors through a ball quotient 
\begin{equation}\label{genus3'}
    p':\   M_3 \rightarrow \bB^6/\Gamma'
\end{equation}
and $p'$ is an open immersion. 

\subsubsection*{g=4}
   A   smooth non-hyperelliptic curve of genus $4$    whose canonical quadric is smooth is a curve $C\subset\mathbb P^1\times\mathbb P^1$   of bidegree $(3,3)$.   For a quadric cone, one uses the resolved cyclic-cover construction of \cite{Kondo02}.   Then a cyclic covering \[ S \rightarrow \bP^1 \times \bP^1 \]
along $C$ will produce a K3 surface $S$ with a order $3$ isomorphism $\tau: S \rightarrow S$. Via the eigen period domain of such K3 surface(see \cite{Kondo02}), there is a holomorphic period map
\begin{equation}\label{genus4}
    p:\ M_g^\circ \rightarrow\   \cD   _6  /\Gamma
\end{equation}
and $p$ factors through 
\begin{equation} \label{genus4'}
    p':\ M_g^\circ \rightarrow\bB^{9}/\Gamma'.
\end{equation}
Let $H_g\subset M_g$ be the locus of hyperelliptic curve of genus $g$. We note that  $H_g$ is a codimension $g-2\ge 2$ subvariety in $M_g$ if $g\ge 4$.  Thus 
\[  \pi_1^{orb}(M_g^\circ)\ \cong \  \pi_1^{orb}(M_g)=\MCG_g \]
where $M_g^\circ:=M_g-H_g$ and $\MCG_g $ is the mapping class group of a compact Riemann surface of genus $g$.

\subsubsection*{g=6}
   Let $M_6^\circ$ denote the locus of nonspecial genus-six curves (neither hyperelliptic, trigonal, bielliptic nor plane quintic). Such a   curve $C$ lies on a unique nodal del Pezzo surface $S$ of degree $5$ as a curve in $|-2K_S|$. By taking the canonical resolution $S'$ of $S$ and a double covering $Y \rightarrow S'$ of $S'$ branched along the curve $C$,  we have a K3 surface $Y$ with    anti-symplectic   involution $\tau: Y \rightarrow Y$, see for example \cite[\S 2 ]{MK11} for more details. Then via the period of K3 surface, there is a    holomorphic   map 
\begin{equation}\label{genus6}
 M_g^\circ \rightarrow\ \cD_T/\Gamma
\end{equation}
where $\cD_T$ is the period domain of the lattice    $T=U^2\oplus E_8(-1)\oplus A_1(-1)^5$ and $\Gamma=O^+(T)$  .

\vspace{0.3cm}

   For the type IV period maps, the differential has complex rank at least two generically, so the argument of Theorem~\ref{h-deform-rigidity_thm_ball_case} proves homotopy deformation rigidity. For the ball-valued factors, a nontrivial homotopy deformation
would give Miyano--Noguchi's totally geodesic product immersion
$E:Z\times W\to Y$ \cite[Theorems~2.15 and 4.9]{MiyanoNoguchi91}.
Both factors have positive dimension: $Z$ parametrizes the deformation,
and the original nonconstant map factors through $W$. A mixed tangent
plane of the product has zero sectional curvature, contradicting the
strictly negative sectional curvature of the ball  .  Then it is natural to ask whether the global rigidity holds for the    period   maps (\ref{genus3}),(\ref{genus3'}),(\ref{genus4}) ,(\ref{genus4'})and (\ref{genus6}) ? 

Note that in the above cases, the classifications of the $\bC$-representation of  $\MCG_g $  due to Korkmaz \cite{korkmaz23}  are only  known for dimension $\le 3g-3$. Then in the view point of group theoretic rigidity,  it will lead to the following interesting question about the  mapping class group 
\begin{question}
   If $g\in \{3,4,6\} $,  is there only one nontrivial group homomorphism $\MCG_g \rightarrow \SO_0(T)$ up to conjugation ?
\end{question}

\subsection{Global Rigidity problems for period maps of  moduli spaces of cubic $n$-folds for  $2\leq n\leq5$} Among hypersurfaces, cubic ones are of particular interest due to their particular geometry, see \cite{Huy}. We denote $\cC_{n}$ the moduli space of smooth cubic $n$-folds for  $2\leq n\leq5$. We make such choice as these cubics have very interesting period maps described as folllows: 
\begin{enumerate}
    \item $n=2$, given $[S]\in \cC_n$,  there is a  natural cyclic covering $Y \rightarrow \bP^3$ of degree $3$ along $S$ and $Y$ is a smooth cubic $3$-fold. Then we have a period map 
    \begin{equation*}
        p:\ \cC_2 \rightarrow \cA_5,\ [S]\ \mapsto [J_Y]
    \end{equation*}
    \item $n=3,5$, there are  period maps via intermediate Jacobian. 
    \item $n=4$,  via primitive cohomology $H^4_{\mathrm{p}}(-,\bZ)$ of a smooth cubic $4$-fold (cf. {\cite[Chapter 6,\S 6]{Huy}}), there is  period map \[   p:\ \cC_4 \rightarrow \cD/\Gamma,\ [X]\ \mapsto [H^{1,3}(X)]  \] where $\cD$ is the period domain associated to the lattice    $T=E_8(-1)^{\oplus2}\oplus U^{\oplus2}\oplus A_2(-1)$ (the negative primitive intersection form) and $\Gamma=\widetilde O^+(T)$  .  
\end{enumerate}
Following the strategy in \S \ref{sec3}, we need to classify the group homomorphism 
\begin{equation}
     \pi_1^{orb}(\cC_3) \rightarrow    \Sp(10,\mathbb Z)
  \end{equation}
up to the conjugation of    $\Sp(10,\mathbb Z)$  .
As the first step to the above problem, we need to understand the 
orbifold fundamental group of the moduli space $M$. Let us recall the Borel construction of  orbifold fundamental group for the orbifold $X/G$ where $X$ is a  complex manifold and $G$ is a reductive group acting on $X$ with finite stabilizer. One may refer to \cite{Looijenga08} for more details. Denote $EG \rightarrow BG$ the universal    principal   $G$-bundle and $X_G:=EG \times_G X$ the $G$-orbit space of $EG \times X$ with $G$  acting diagonally. Then $\pi_1^{orb}(X/G)=\pi_1(X_G)$. From the construction there is a topological    fibration   $ X\rightarrow  X_G  \rightarrow BG$, which induces Serre's long exact sequence
\begin{equation*}
    \cdots   \rightarrow   \pi_1(G) \rightarrow  \pi_1(X) \rightarrow    \pi_1^{orb}(X/G) \rightarrow \pi_0(G) \rightarrow 1.
\end{equation*}
Now we apply the above formulation to $G=\PGL(n+1,\bC)$ and $ U\subset |L|$ the locus of smooth varieties.  As $G=\PGL(n+1,\bC)$ is connected,  we have a surjective morphism $$ \pi_1(U)\twoheadrightarrow \pi_1^{orb}(U/G)$$
where the Kernel is given by $\im( \pi_1(\PGL(n+1,\bC)) \rightarrow  \pi_1(U,p))$ and it is induced by the map \[ \PGL(n+1,\bC) \rightarrow U,\ \ g\mapsto gp   .\]
By \cite[Main theorem]{Lonne09}, $ \pi_1(U)$ is a generalization of Braid group, where the generator of $ \pi_1(U)$ is known. We ask the following 
natural question:

\begin{question}
    what is the orbifold fundamental group $\pi_1^{orb}(M)$ ?  
\end{question}

Now we consider the orbifold obstruction to global rigidity for period maps of these cubic $n$-folds. Note that in \cite{BFJS26} Baldi-Farb-Javanpeykar-Stover proved that $\aut_{orb}(\cC_2)=\{1\}$, suggesting that the following question may have a positive answer.  

\begin{question}\label{GRC3}
    The non-constant holomorphic map from $\cC_2$ or  $\cC_3$ to $\cA_5$ must be the period map ? 
\end{question}
\begin{rem}
   The presentation of the orbifold fundamental group of moduli space of smooth cubic surfaces is known by \cite{Looijenga08} and also recent work \cite{ABL26}.  Their methods may be useful to understand the above question. We leave it for further study. 
\end{rem}
Indeed  as $\cC_2= \bB^{4}/\Gamma-H$ for a    Heegner   divisor $H$ (see \cite[Theorem 2.20]{ACT2002}),  Baldi-Farb-Javanpeykar-Stover \cite{BFJS26} reduced  $\aut_{orb}(\cC_2)=\{1\}$ to show $\aut_{orb}(\bB^{4}/\Gamma )=\{\mathrm{Id}\}$ by studying   $\pi_1^{orb}(\bB^{4}/\Gamma-H)$. To give more positive evidence for  Question \ref{GRC3} when the source moduli space is  $\cC_3$,  we give a generalization to  Eisenstein–Picard modular $n$-folds.  Before procedure we recall some basic definitions. Let $\cE$ be the  Eisenstein ring, i.e., the ring of integer of the imaginary quadratic field $\bQ(\xi)=\bQ(\sqrt{-3})$ where $\xi=\exp(\frac{2\pi\sqrt{-1}}{3})$. Denote $(\cE^{n+1},h)$ be the Hermitian space where the Hermitian form is given by \[h(z)=-\overline{z}_0z_0+\overline{z}_1z_1+\cdots +\overline{z}_nz_n.\]
Let $\pi: \mathrm{SU}(n,1)  \rightarrow \mathrm{PU}(n,1) $ be the natural projection map and \[\mathrm{U}(\cE,h):=\{ A\in \GL(\cE^{n+1})\ |\ h(-,-)=h(A-,A-) \}\] and  $ \mathrm{SU}(\cE,h):=\{A\in \mathrm{U}(\cE,h)\ |\ \det(A)=1\} $,  then \[\Gamma_n:=\pi(\mathrm{SU}(\cE,h)) \] is an arithmetic subgroup in $\mathrm{PU}(n,1)$, which acts on the ball \[\bB^n:=\{z\in \bP^n\ |\ h|_z   <  0 \}.\]
The quotient space $\bB^n/\Gamma_n$ is a quasi-projective variety of  dimension $n$ with orbifold singularities, which is called an Eisenstein–Picard modular $n$-fold.

\begin{pro}\label{Eisenstein–Picard}
     For  $n\geq 2$, the orbifold automorphism group of    the   Eisenstein–Picard modular $n$-fold   above   is given by
\begin{equation*}
    \aut_{orb}({\bB^{n}/\Gamma_n})=C_{g.c.d(6,n+1)}
\end{equation*}
where $C_{g.c.d(6,n+1)}$ is a cyclic group of order $g.c.d(6,n+1)$. 
\end{pro}

\begin{proof}
First we claim that the normalizer 
\begin{equation*}
    N_{\mathrm{PU}(n,1)}(\Gamma_n)=\mathrm{PU}(\cE,h).
\end{equation*}
Assume the claim, then 
\[ \aut_{orb}({\bB^{n}/\Gamma_n})=\frac{N_{\mathrm{PU}(n,1)}(\Gamma_n)}{\Gamma_n}=\frac{\mathrm{PU}(\cE,h)}{\Gamma_n} \]
By the projection $\pi: \mathrm{SU}(n,1)  \rightarrow \mathrm{PU}(n,1) $,  \[ \frac{\mathrm{PU}(\cE,h)}{\Gamma_n}\cong \frac{\mu_6}{\mu_6^{n+1}}\]
where $\mu_6^{n+1} \le \mu_6$ is the subgroup.
Then the proposition follows from the fact that 
$\frac{\mu_6}{\mu_6^{n+1}} \cong C_{g.c.d(6,n+1)}$.

Now we divide the proof of the claim into following two steps.

\paragraph{Step1:} It is sufficient to show that 
\begin{equation} \label{KeySU}
   N_{\mathrm{SU}(n,1)}(\pi^{-1}(\Gamma_n))=\pi^{-1}(\mathrm{PU}(\cE,h)).
\end{equation}

Indeed as  $ N_{\mathrm{SU}(n,1)}(\pi^{-1}(\Gamma_n)) \le \mathrm{SU}(n,1)$ is a maximal arithmetic lattice by \cite[Proposition 3.8]{ES14} and $\pi$ is surjective, then by equality (\ref{KeySU})  $\pi(N_{\mathrm{SU}(n,1)}(\widetilde{\Gamma}_n))= \mathrm{PU}(\cE,h) $ is also a maximal arithmetic lattice in $\mathrm{PU}(n,1)$. One can easily check that $\mathrm{PU}(\cE,h) \subset N_{\mathrm{PU}(n,1)}(\Gamma_n)$ with finite index and thus  $N_{\mathrm{PU}(n,1)}(\Gamma_n)$ is also an arithmetic lattice, then the maximality  of $\mathrm{PU}(\cE,h)$ forces $\mathrm{PU}(\cE,h) = N_{\mathrm{PU}(n,1)}(\Gamma_n)$.

\paragraph{Step2:} Note that \[\pi^{-1}(\mathrm{PU}(\cE,h)) =\{ \lambda A\ |\ \lambda^{n+1} \cdot \det(A)=1, \ \ A\in \mathrm{U}(\cE,h)\} \] and $\pi^{-1}(\mathrm{PU}(\cE,h)) \subset N_{\mathrm{SU}(n,1)}(\pi^{-1}(\Gamma_n))$,  then there is an exact sequence 
\begin{equation}\label{SES:1}
    1 \rightarrow \frac{\ker(\pi)}{\ker(\pi)\cap \mathrm{SU}(\cE,h)}\rightarrow    \frac{\pi^{-1}(\mathrm{PU}(\cE,h)}{\mathrm{SU}(\cE,h)}  \rightarrow \frac{\mathrm{PU}(\cE,h)}{\Gamma_n}  \rightarrow 1.
\end{equation}
By direct computation, 
\[\ker(\pi)\cap \mathrm{SU}(\cE,h)=\{ \ \lambda \mathrm{Id}\ |\   \lambda\in \bQ(\sqrt{-3}),\ \lambda^{n+1}=1 \}\cong \mu_{n+1} \cap \mu_6.\]
Thus the exact sequence (\ref{SES:1}) shows the index upper  bound
\begin{equation*}
    \begin{split}
        [ N_{\mathrm{SU}(n,1)}(\pi^{-1}(\Gamma_n)):\mathrm{SU}(\cE,h)]  &\ge  [\pi^{-1}(\mathrm{PU}(\cE,h): \mathrm{SU}(\cE,h)]\\
        &=|\frac{\mu_{n+1}}{\mu_{n+1}\cap \mu_6}|\cdot |\frac{\mu_6}{\mu_6^{n+1}}|  = n+1.
    \end{split}
\end{equation*}

Then to show the equality (\ref{KeySU}), it remains to  show the index lower bound $ [ N_{\mathrm{SU}(n,1)}(\pi^{-1}(\Gamma_n)):\mathrm{SU}(\cE,h)]\le n+1$. This can be done by following Emery–Stover's index computations for normalizers of principal arithmetic lattices in \cite[\S 3.2]{ES14}:
\begin{itemize}
    \item if $n$ is even, then $\mathrm{SU}(\cE,h)$ is quasi-split at every finite prime and thus \cite[Lemma 3.7]{ES14} implies \[  N_{\mathrm{SU}(n,1)}(\pi^{-1}(\Gamma_n)):\mathrm{SU}(\cE,h)]\le  (n+1)h_{\bQ(\sqrt{-3})}=n+1\] 
as  the class number $h_{\bQ(\sqrt{-3})}$ of the number field $\bQ(\sqrt{-3})$ is  $1$.
    \item if $n$ is odd, then $\mathrm{SU}(\cE,h)$ is  not quasi-split at the set $T=\{3\}$ but is ramified over $T=\{3\}$, then \cite[Lemma 3.7]{ES14} can also be applied directly, which also implies 
    $ [N_{\mathrm{SU}(n,1)}(\pi^{-1}(\Gamma_n)):\mathrm{SU}(\cE,h)]\le n+1$.
\end{itemize}

\end{proof}
In particular, the Ball quotient containing the period of intermediate Jacobian of smooth  cubic $3$-folds has only trivial orbifold automorphism. 



\subsection*{Rigidity of period maps via Albanese variety}
Suppose    $\pi:X\to M$ is a smooth projective family with a relative polarization   over  $M$ such that $M$ has orbifold singularities at worst and $g:=\dim H^0(X_t,\Omega_{X_t})>0$ for any $t\in M$.   The relative polarization induces a fixed polarization type $\delta$ on the Albanese varieties over each connected component of $M$.   Then there is a period map \begin{equation} \label{periodmap:alb}
     p:\ M   \to\mathcal A_{g,\delta}  ,   \quad   t   \mapsto[\operatorname{Alb}(X_t)].
\end{equation}    The target is $\mathcal A_g$ when the induced polarization is principal.
Here   $$\Alb(X_t):=\frac{ H^0(X_t,\Omega_{X_t})^\vee}{ H_1(X_t,\bZ)}$$ is the Albanese variety of $X_t$. This can be viewed as a direct generalization of the Torelli map as period map of smooth projective curves. 
we may ask 
\begin{question}
    when the period map in (\ref{periodmap:alb}) is globally rigid ?
\end{question}

  \vspace{0.3cm}
\bibliographystyle{alpha}
   \bibliography{period_maps_minimal_20260926}

@article{BorelNarasimhan66,
  author  = {Borel, Armand and Narasimhan, Raghavan},
  title   = {Uniqueness Conditions for Certain Holomorphic Mappings},
  journal = {Inventiones Mathematicae},
  volume  = {2},
  year    = {1966},
  pages   = {247--255}
}

@article{Borel72,
  author  = {Borel, Armand},
  title   = {Some Metric Properties of Arithmetic Quotients of Symmetric Spaces and an Extension Theorem},
  journal = {Journal of Differential Geometry},
  volume  = {6},
  year    = {1972},
  pages   = {543--560}
}

@misc{Peters00,
  author        = {Peters, Chris A. M.},
  title         = {Arakelov-Type Inequalities for Hodge Bundles},
  year          = {2000},
  eprint        = {math/0007102},
  archivePrefix = {arXiv},
  primaryClass  = {math.AG}
}

@article {Ric88,
    AUTHOR = {Richardson, R. W.},
     TITLE = {Conjugacy classes of {$n$}-tuples in {L}ie algebras and
              algebraic groups},
   JOURNAL = {Duke Math. J.},
  FJOURNAL = {Duke Mathematical Journal},
    VOLUME = {57},
      YEAR = {1988},
    NUMBER = {1},
     PAGES = {1--35},
      ISSN = {0012-7094,1547-7398},
   MRCLASS = {20G15 (14L30 17B45 22E46)},
  MRNUMBER = {952224},
MRREVIEWER = {Michel\ Brion},
       DOI = {10.1215/S0012-7094-88-05701-8},
       URL = {https://doi.org/10.1215/S0012-7094-88-05701-8},
}

@article{Farb22,
  title={Global rigidity of the period mapping},
  author={Farb, Benson},
  journal={Journal of Algebraic Geometry},
  year={2022}
}

@incollection {Farb24,
    AUTHOR = {Farb, Benson},
     TITLE = {Rigidity of moduli spaces and algebro-geometric constructions},
 BOOKTITLE = {Surveys in differential geometry 2021. {C}hern: a great
              geometer of the 20th century},
    SERIES = {Surv. Differ. Geom.},
    VOLUME = {26},
     PAGES = {31--49},
 PUBLISHER = {Int. Press, Boston, MA},
      YEAR = {2024},
      ISBN = {978-1-57146-464-4},
   MRCLASS = {14D20 (58D27)},
  MRNUMBER = {4748488},
}

@incollection {Kondo02,
    AUTHOR = {Kondo, Shigeyuki},
     TITLE = {The moduli space of curves of genus 4 and {D}eligne-{M}ostow's
              complex reflection groups},
 BOOKTITLE = {Algebraic geometry 2000, {A}zumino ({H}otaka)},
    SERIES = {Adv. Stud. Pure Math.},
    VOLUME = {36},
     PAGES = {383--400},
 PUBLISHER = {Math. Soc. Japan, Tokyo},
      YEAR = {2002},
      ISBN = {4-931469-20-5},
   MRCLASS = {14H15 (14D07 14H45 14J28 32S40 33C80)},
  MRNUMBER = {1971521},
MRREVIEWER = {I.\ Dolgachev},
       DOI = {10.2969/aspm/03610383},
       URL = {https://doi.org/10.2969/aspm/03610383},
}

@article {MK11,
    AUTHOR = {Artebani, Michela and Kondo, Shigeyuki},
     TITLE = {The moduli of curves of genus six and {$K3$} surfaces},
   JOURNAL = {Trans. Amer. Math. Soc.},
  FJOURNAL = {Transactions of the American Mathematical Society},
    VOLUME = {363},
      YEAR = {2011},
    NUMBER = {3},
     PAGES = {1445--1462},
      ISSN = {0002-9947,1088-6850},
   MRCLASS = {14H10 (14J28 14M27)},
  MRNUMBER = {2737272},
MRREVIEWER = {Lars\ Halvard\ Halle},
       DOI = {10.1090/S0002-9947-2010-05126-8},
       URL = {https://doi.org/10.1090/S0002-9947-2010-05126-8},
}

@article {Michela09,
    AUTHOR = {Artebani, Michela},
     TITLE = {A compactification of $M_3$ via $K3$ surfaces},
   JOURNAL = {Nagoya Math. J.},
  FJOURNAL = {Nagoya Mathematical Journal},
    VOLUME = {196},
      YEAR = {2009},
     PAGES = {1--26},
      ISSN = {0027-7630,2152-6842},
   MRCLASS = {14J10 (14H10 14J28)},
  MRNUMBER = {2591089},
MRREVIEWER = {Eleonora\ Palmieri},
       DOI = {10.1017/S0027763000009776},
       URL = {https://doi.org/10.1017/S0027763000009776},
}

@article {Kondo00,
    AUTHOR = {Kondo, Shigeyuki},
     TITLE = {A complex hyperbolic structure for the moduli space of curves
              of genus three},
   JOURNAL = {J. Reine Angew. Math.},
  FJOURNAL = {Journal f\"{u}r die Reine und Angewandte Mathematik. [Crelle's
              Journal]},
    VOLUME = {525},
      YEAR = {2000},
     PAGES = {219--232},
      ISSN = {0075-4102,1435-5345},
   MRCLASS = {14H15 (14J28)},
  MRNUMBER = {1780433},
MRREVIEWER = {Dan\ Avritzer},
       DOI = {10.1515/crll.2000.069},
       URL = {https://doi.org/10.1515/crll.2000.069},
}

@article {Nog88,
    AUTHOR = {Noguchi, Junjiro},
     TITLE = {Moduli spaces of holomorphic mappings into hyperbolically
              imbedded complex spaces and locally symmetric spaces},
   JOURNAL = {Invent. Math.},
  FJOURNAL = {Inventiones Mathematicae},
    VOLUME = {93},
      YEAR = {1988},
    NUMBER = {1},
     PAGES = {15--34},
      ISSN = {0020-9910,1432-1297},
   MRCLASS = {32H20 (32G13 32H25 32M15)},
  MRNUMBER = {943922},
MRREVIEWER = {R.\ R.\ Simha},
       DOI = {10.1007/BF01393686},
       URL = {https://doi.org/10.1007/BF01393686},
}

@article {SchoenYau79,
    AUTHOR = {Schoen, Richard and Yau, Shing Tung},
     TITLE = {Compact group actions and the topology of manifolds with
              nonpositive curvature},
   JOURNAL = {Topology},
  FJOURNAL = {Topology. An International Journal of Mathematics},
    VOLUME = {18},
      YEAR = {1979},
    NUMBER = {4},
     PAGES = {361--380},
      ISSN = {0040-9383},
   MRCLASS = {53C20 (58E20)},
  MRNUMBER = {551017},
MRREVIEWER = {Helmut\ Kaul},
       DOI = {10.1016/0040-9383(79)90026-0},
       URL = {https://doi.org/10.1016/0040-9383(79)90026-0},
}

@article {Saito-Zucker,
    AUTHOR = {Saito, Masa-Hiko and Zucker, Steven},
     TITLE = {Classification of nonrigid families of {$K3$} surfaces and a
              finiteness theorem of {A}rakelov type},
   JOURNAL = {Math. Ann.},
  FJOURNAL = {Mathematische Annalen},
    VOLUME = {289},
      YEAR = {1991},
    NUMBER = {1},
     PAGES = {1--31},
      ISSN = {0025-5831,1432-1807},
   MRCLASS = {14J10 (14G40 14J28 32J25)},
  MRNUMBER = {1087233},
MRREVIEWER = {Gerd\ Faltings},
       DOI = {10.1007/BF01446555},
       URL = {https://doi.org/10.1007/BF01446555},
}

@book {Mok89,
    AUTHOR = {Mok, Ngaiming},
     TITLE = {Metric rigidity theorems on {H}ermitian locally symmetric
              manifolds},
    SERIES = {Series in Pure Mathematics},
    VOLUME = {6},
 PUBLISHER = {World Scientific Publishing Co., Inc., Teaneck, NJ},
      YEAR = {1989},
     PAGES = {xiv+278},
      ISBN = {9971-50-800-1; 9971-50-802-8},
   MRCLASS = {32M15 (32C17 32H20 32L07 53C35)},
  MRNUMBER = {1081948},
MRREVIEWER = {B.\ Gilligan},
       DOI = {10.1142/0773},
       URL = {https://doi.org/10.1142/0773},
}

@article {Sankaran97,
    AUTHOR = {Sankaran, G. K.},
     TITLE = {Moduli of polarised abelian surfaces},
   JOURNAL = {Math. Nachr.},
  FJOURNAL = {Mathematische Nachrichten},
    VOLUME = {188},
      YEAR = {1997},
     PAGES = {321--340},
      ISSN = {0025-584X,1522-2616},
   MRCLASS = {14K10 (11F46)},
  MRNUMBER = {1484680},
MRREVIEWER = {Jae-Hyun\ Yang},
       DOI = {10.1002/mana.19971880117},
       URL = {https://doi.org/10.1002/mana.19971880117},
}

@incollection {Milne,
    AUTHOR = {Milne, J. S.},
     TITLE = {Introduction to {S}himura varieties},
 BOOKTITLE = {Harmonic analysis, the trace formula, and {S}himura varieties},
    SERIES = {Clay Math. Proc.},
    VOLUME = {4},
     PAGES = {265--378},
 PUBLISHER = {Amer. Math. Soc., Providence, RI},
      YEAR = {2005},
      ISBN = {0-8218-3844-X},
   MRCLASS = {11G18 (14G35)},
  MRNUMBER = {2192012},
MRREVIEWER = {Ulrich\ G\"ortz},
}

@book{Huybrechts2016,
  title={Lectures on K3 surfaces},
  author={Huybrechts, Daniel},
  volume={158},
  year={2016},
  publisher={Cambridge University Press}
}

@article {BHC62,
    AUTHOR = {Borel, Armand and Harish-Chandra},
     TITLE = {Arithmetic subgroups of algebraic groups},
   JOURNAL = {Ann. of Math. (2)},
  FJOURNAL = {Annals of Mathematics. Second Series},
    VOLUME = {75},
      YEAR = {1962},
     PAGES = {485--535},
       DOI = {10.2307/1970210},
       URL = {https://doi.org/10.2307/1970210},
}

@incollection {Deligne87,
    AUTHOR = {Deligne, Pierre},
     TITLE = {Un th\'{e}or\`{e}me de finitude pour la monodromie},
 BOOKTITLE = {Discrete groups in geometry and analysis (New Haven, Conn., 1984)},
    EDITOR = {Howe, Roger},
    SERIES = {Progr. Math.},
    VOLUME = {67},
     PAGES = {1--19},
 PUBLISHER = {Birkh\"auser Boston, Boston, MA},
      YEAR = {1987},
       DOI = {10.1007/978-1-4899-6664-3_1},
       URL = {https://doi.org/10.1007/978-1-4899-6664-3_1},
}

@book{ACT11,
  title={The moduli space of cubic threefolds as a ball quotient},
  author={Allcock, Daniel and Carlson, James and Toledo, Domingo},
  series={Memoirs of the American Mathematical Society},
  volume={209, no. 985},
  year={2011},
  publisher={American Mathematical Society}
}

@phdthesis{Reid72,
  title={The complete intersection of two or more quadrics},
  author={Reid, Miles A},
  year={1972},
  school={University of Cambridge Cambridge}
}

@article{Collino86,
  title={The Abel-Jacobi isomorphism for the cubic fivefold},
  author={Collino, Alberto},
  journal={Pacific Journal of Mathematics},
  volume={122},
  number={1},
  pages={43--55},
  year={1986},
  publisher={Mathematical Sciences Publishers}
}

@article {CL18,
    AUTHOR = {Cr\'{e}tois, R\'{e}mi and Lang, Lionel},
     TITLE = {The vanishing cycles of curves in toric surfaces {I}},
   JOURNAL = {DukeCompos. Math. J.},
  FJOURNAL = {DukeCompositio Mathematical Journal},
    VOLUME = {57154},
      YEAR = {1982018},
    NUMBER = {18},
     PAGES = {1--35},
       DOI = {10.1215/S0012-7094-88-05701-8659--1697},
      ISSN = {0010-437X,1570-5846},
   MRCLASS = {14M25 (14T05 32S30)},
  MRNUMBER = {3830549},
MRREVIEWER = {Frank\ Sottile},
       DOI = {10.1112/s0010437x18007200},
       URL = {https://doi.org/10.1215/S0012-7094-88-05701-8112/s0010437x18007200},
}

@article{FM98,
  title={Superrigidity and mapping class groups},
  author={Farb, Benson and Masur, Howard},
  journal={Topology},
  volume={37},
  number={6},
  pages={1169--1176},
  year={1998},
  publisher={Elsevier}
}

@article {BDS11,
    AUTHOR = {Behrstock, Jason and Dru\c{t}u, Cornelia and Sapir, Mark},
     TITLE = {Addendum: {M}edian structures on asymptotic cones and
              homomorphisms into mapping class groups [MR2783135]},
   JOURNAL = {Proc. Lond. Math. Soc. (3)},
  FJOURNAL = {Proceedings of the London Mathematical Society. Third Series},
    VOLUME = {95102},
      YEAR = {200211},
    NUMBER = {3},
     PAGES = {19555--5862},
}

@article{Looijenga08,
  title={Artin groups and the fundamental groups of some moduli spaces},
  author={Looijenga, Eduard},
  journal={Journal of Topology},
  volume={1},
  number={1},
  pages={187--216},
  year={2008},
  publisher={Wiley Online Library},
}

@article {Lonne09,
    AUTHOR = {L\"{o}nne, Michael},
     TITLE = {Fundamental groups of projective discriminant complements},
   JOURNAL = {Duke Math. J.},
  FJOURNAL = {Duke Mathematical Journal},
    VOLUME = {150},
      YEAR = {2009},
    NUMBER = {2},
     PAGES = {357--405},
      ISSN = {0012-7094,1547-7398},
   MRCLASS = {14F35 (14D05 14J70)},
  MRNUMBER = {2569617},
MRREVIEWER = {Michael\ Friedman},
       DOI = {10.1215/00127094-2009-055},
       URL = {https://doi.org/10.1215/00127094-2009-055},
}

@article {KS22,
    AUTHOR = {Krieg, Aloys and Schaps, Felix},
     TITLE = {Maximal discrete subgroups of {$SO^+(2,n+2)$}},
   JOURNAL = {Proc. Amer. Math. Soc.},
  FJOURNAL = {Proceedings of the American Mathematical Society},
    VOLUME = {150},
      YEAR = {2022},
    NUMBER = {6},
     PAGES = {2357--2365},
      ISSN = {0002-9939,1088-6826},
   MRCLASS = {11F06 (11F55 22E40)},
  MRNUMBER = {4399255},
MRREVIEWER = {Haowu\ Wang},
       DOI = {10.1090/proc/15889},
       URL = {https://doi.org/10.1090/proc/15889},
}

@article{BFJS26,
  title={Automorphisms of the moduli space of smooth cubic surfaces and its fundamental group},
  author={Baldi, Gregorio and Farb, Benson and Javanpeykar, Ariyan and Stover, Matthew},
  journal={arXiv preprint arXiv:2605.16658},
  year={2026}
}

@book {Huy,
    AUTHOR = {Huybrechts, Daniel},
     TITLE = {The geometry of cubic hypersurfaces},
    SERIES = {Cambridge Studies in Advanced Mathematics},
    VOLUME = {206},
 PUBLISHER = {Cambridge University Press, Cambridge},
      YEAR = {2023},
     PAGES = {xvii+441},
      ISBN = {978-1-009-28000-6; [9781009280020]},
   MRCLASS = {14J70 (14C30 14J30 14J35)},
  MRNUMBER = {4589520},
MRREVIEWER = {Fei\ Hu},
}

@article{Niku79,
  author = {Nikulin, V. V.},
  title = {Integral symmetric bilinear forms and some of their applications},
  journal = {Math. USSR-Izv.},
  volume = {14},
  number = {1},
  year = {1980},
  pages = {103--167},
  doi = {10.1070/IM1980v014n01ABEH001060}
}

@article {ES14,
    AUTHOR = {Emery, Vincent and Stover, Matthew},
     TITLE = {Covolumes of nonuniform lattices in {${\rm PU}(n,1)$}},
   JOURNAL = {Amer. J. Math.},
  FJOURNAL = {American Journal of Mathematics},
    VOLUME = {136},
      YEAR = {2014},
    NUMBER = {1},
     PAGES = {143--164},
      ISSN = {0002-9327,1080-6377},
   MRCLASS = {22E40 (11E57 14L30 51M10 53C23)},
  MRNUMBER = {3163356},
MRREVIEWER = {Martin\ Deraux},
       DOI = {10.1353/ajm.2014.0007},
       URL = {https://doi.org/10.1353/ajm.2014.0007},
}

@article {Ihara67,
    AUTHOR = {Ihara, Shin-ichiro},
     TITLE = {Holomorphic imbeddings of symmetric domains},
   JOURNAL = {J. Math. Soc. Japan},
  FJOURNAL = {Journal of the Mathematical Society of Japan},
    VOLUME = {19},
      YEAR = {1967},
     PAGES = {261--302},
      ISSN = {0025-5645,1881-1167},
   MRCLASS = {32.32 (22.00)},
  MRNUMBER = {214807},
MRREVIEWER = {Ichiro\ Satake},
       DOI = {10.2969/jmsj/01930261},
       URL = {https://doi.org/10.2969/jmsj/01930261},
}

@article {Satake65,
    AUTHOR = {Satake, Ichir\^o},
     TITLE = {Holomorphic imbeddings of symmetric domains into a {S}iegel
              space},
   JOURNAL = {Amer. J. Math.},
  FJOURNAL = {American Journal of Mathematics},
    VOLUME = {87},
      YEAR = {1965},
     PAGES = {425--461},
      ISSN = {0002-9327,1080-6377},
   MRCLASS = {32.65 (10.23)},
  MRNUMBER = {196134},
MRREVIEWER = {A.\ Kor\'anyi},
       DOI = {10.2307/2373012},
       URL = {https://doi.org/10.2307/2373012},
}

@incollection {Beauville,
    AUTHOR = {Beauville, Arnaud},
     TITLE = {Le groupe de monodromie des familles universelles
              d'hypersurfaces et d'intersections compl\`etes},
 BOOKTITLE = {Complex analysis and algebraic geometry ({G}\"{o}ttingen,
              1985)},
    SERIES = {Lecture Notes in Math.},
    VOLUME = {1194},
     PAGES = {8--18},
 PUBLISHER = {Springer, Berlin},
      YEAR = {1986},
      ISBN = {3-540-16490-1},
   MRCLASS = {14J10 (14D05 14H10 14M10)},
  MRNUMBER = {855873},
MRREVIEWER = {I.\ Dolgachev},
       DOI = {10.1007/BFb0076991},
       URL = {https://doi.org/10.1007/BFb0076991},
}

@article {Deraux11,
    AUTHOR = {Deraux, Martin},
     TITLE = {Forgetful maps between {D}eligne-{M}ostow ball quotients},
   JOURNAL = {Geom. Dedicata},
  FJOURNAL = {Geometriae Dedicata},
    VOLUME = {150},
      YEAR = {2011},
     PAGES = {377--389},
      ISSN = {0046-5755,1572-9168},
   MRCLASS = {32Q05 (14D06 22E40 32J99)},
  MRNUMBER = {2753711},
MRREVIEWER = {G.\ K.\ Sankaran},
       DOI = {10.1007/s10711-010-9511-x},
       URL = {https://doi.org/10.1007/s10711-010-9511-x},
}

@article{korkmaz23,
  title={Low-dimensional linear representations of mapping class groups},
  author={Korkmaz, Mustafa},
  journal={Journal of Topology},
  volume={16},
  number={3},
  pages={899--935},
  year={2023},
  publisher={Wiley Online Library}
}

@article{servan25,
  title={On the uniqueness of the Prym map},
  author={Serv{\'a}n, Carlos A},
  journal={Journal of Differential Geometry},
  volume={131},
  number={3},
  pages={801--830},
  year={2025},
  publisher={Lehigh University}
}

@article{DK06,
  title={Moduli of K3 surfaces and complex ball quotients},
  author={Dolgachev, Igor V and Kondo, Shigeyuki},
  journal={PROGRESS IN MATHEMATICS-BOSTON-},
  volume={260},
  pages={43},
  year={2006},
  publisher={Birkh{\"a}user}
}

@article {BB66,
    AUTHOR = {Baily, Jr., W. L. and Borel, A.},
     TITLE = {Compactification of arithmetic quotients of bounded symmetric
              domains},
   JOURNAL = {Ann. of Math. (2)},
  FJOURNAL = {Annals of Mathematics. Second Series},
    VOLUME = {84},
      YEAR = {1966},
     PAGES = {442--528},
      ISSN = {0003-486X},
   MRCLASS = {32.65},
  MRNUMBER = {216035},
MRREVIEWER = {A.\ Kor\'{a}nyi},
       DOI = {10.2307/1970457},
       URL = {https://doi.org/10.2307/1970457},
}

@book {DF04,
    AUTHOR = {Dummit, David S. and Foote, Richard M.},
     TITLE = {Abstract algebra},
   EDITION = {Third},
 PUBLISHER = {John Wiley \& Sons, Inc., Hoboken, NJ},
      YEAR = {2004},
     PAGES = {xii+932},
      ISBN = {0-471-43334-9},
   MRCLASS = {00-01 (16-01 20-01)},
  MRNUMBER = {2286236},
}

@ARTICLE{Zhong26,
       author = {{Zhong}, Xiyan},
        title = "{Rigidity of the period map up to finite covers}",
      journal = {arXiv e-prints},
         year = 2026,
        month = aug,
          eid = {arXiv:2608.29351},
        pages = {arXiv:2608.29351},
archivePrefix = {arXiv},
       eprint = {2608.29351},
 primaryClass = {math.GT},
       adsurl = {https://ui.adsabs.harvard.edu/abs/2026arXiv260829351Z}
}

@book {Gru07,
    AUTHOR = {Gruber, Peter M.},
     TITLE = {Convex and discrete geometry},
    SERIES = {Grundlehren der mathematischen Wissenschaften [Fundamental
              Principles of Mathematical Sciences]},
    VOLUME = {336},
 PUBLISHER = {Springer, Berlin},
      YEAR = {2007},
     PAGES = {xiv+578},
      ISBN = {978-3-540-71132-2},
   MRCLASS = {52-02 (11-02 46B20 49-02 49J53 90C25)},
  MRNUMBER = {2335496},
MRREVIEWER = {Aleksandr\ Koldobsky},
}

@article{JostYau93,
  author  = {Jost, J{\"u}rgen and Yau, Shing-Tung},
  title   = {Harmonic mappings and algebraic varieties over function fields},
  journal = {American Journal of Mathematics},
  volume  = {115},
  number  = {6},
  year    = {1993},
  pages   = {1197--1227},
  doi     = {10.2307/2374964}
}

@article{Par68,
  author  = {Parshin, A. N.},
  title   = {Algebraic curves over function fields. I},
  journal = {Mathematics of the USSR-Izvestiya},
  volume  = {2},
  number  = {5},
  year    = {1968},
  pages   = {1145--1170},
  doi     = {10.1070/IM1968v002n05ABEH000723}
}

@article{Ara71,
  author  = {Arakelov, S. Yu.},
  title   = {Families of algebraic curves with fixed degeneracies},
  journal = {Mathematics of the USSR-Izvestiya},
  volume  = {5},
  number  = {6},
  year    = {1971},
  pages   = {1277--1302},
  doi     = {10.1070/IM1971v005n06ABEH001235}
}

@article{Bender80,
  author = {Bender, Peter},
  title = {{Eine Pr{\"a}sentation der symplektischen Gruppe Sp(4, Z) mit 2 Erzeugenden und 8 definierenden Relationen}},
  journal = {J. Algebra},
  volume = {65},
  number = {2},
  year = {1980},
  pages = {328--331}
}

@incollection{Mok12,
  author = {Mok, Ngaiming},
  title = {Projective algebraicity of minimal compactifications of
           complex-hyperbolic space forms of finite volume},
  booktitle = {Perspectives in Analysis, Geometry, and Topology},
  series = {Progress in Mathematics},
  volume = {296},
  pages = {331--354},
  publisher = {Birkh{\"a}user/Springer},
  address = {New York},
  year = {2012},
  doi = {10.1007/978-0-8176-8277-4_14}
}

@incollection{MiyanoNoguchi91,
  AUTHOR = {Miyano, Toshiki and Noguchi, Junjiro},
  TITLE = {Moduli spaces of harmonic and holomorphic mappings and
           {Diophantine} geometry},
  BOOKTITLE = {Prospects in Complex Geometry ({Katata} and {Kyoto}, 1989)},
  EDITOR = {Noguchi, Junjiro and Ohsawa, Takeo},
  SERIES = {Lecture Notes in Mathematics},
  VOLUME = {1468},
  PAGES = {227--253},
  PUBLISHER = {Springer},
  ADDRESS = {Berlin},
  YEAR = {1991},
  DOI = {10.1007/BFb0086196}
}

@article{To89,
  AUTHOR = {To, Wing-Keung},
  TITLE = {Hermitian metrics of semi-negative curvature on quotients
           of bounded symmetric domains},
  JOURNAL = {Invent. Math.},
  VOLUME = {95},
  NUMBER = {3},
  YEAR = {1989},
  PAGES = {559--578},
  DOI = {10.1007/BF01393890}
}

@article{DHNT15,
  title={Families of lattice polarized K3 surfaces with monodromy},
  author={Doran, Charles F and Harder, Andrew and Novoseltsev, Andrey Y and Thompson, Alan},
  journal={International mathematics research notices},
  volume={2015},
  number={23},
  pages={12265--12318},
  year={2015},
  publisher={Oxford University Press}
}

@article {ACT2002,
    AUTHOR = {Allcock, Daniel and Carlson, James A. and Toledo, Domingo},
     TITLE = {The complex hyperbolic geometry of the moduli space of cubic
              surfaces},
   JOURNAL = {J. Algebraic Geom.},
  FJOURNAL = {Journal of Algebraic Geometry},
    VOLUME = {11},
      YEAR = {2002},
    NUMBER = {4},
     PAGES = {659--724},
      ISSN = {1056-3911,1534-7486},
   MRCLASS = {32G20 (11F23 14J10 14J15 14K30 32M15 32Q45)},
  MRNUMBER = {1910264},
MRREVIEWER = {I.\ Dolgachev},
       DOI = {10.1090/S1056-3911-02-00314-4},
       URL = {https://doi.org/10.1090/S1056-3911-02-00314-4},
}

@article{ABL26,
  title={The moduli space of cubic surfaces and its orbifold fundamental group},
  author={Allcock, Daniel and Basak, Tathagata and Looijenga, Eduard},
  journal={arXiv preprint arXiv:2610.16080},
  year={2026}
}
  \end{document}